\documentclass{article}
\usepackage[top=1in,bottom=1in,left=1.3in,right=1.4in]{geometry}
\usepackage[utf8]{inputenc} % allow utf-8 input
\usepackage[T1]{fontenc}    % use 8-bit T1 fonts
\usepackage{hyperref}       % hyperlinks
\usepackage{url}            % simple URL typesetting
\usepackage{booktabs}       % professional-quality tables
\usepackage{amsfonts}       % blackboard math symbols
\usepackage{nicefrac}       % compact symbols for 1/2, etc.
\usepackage{xcolor}         % colors
\usepackage{amsmath, amssymb,amsfonts,bbm,verbatim,multirow,color}
\usepackage{amsthm}
\usepackage{amscd}

\usepackage{comment}

\newtheorem{thm}{Theorem}[section]

\newtheorem{prp}[thm]{Proposition}
\newtheorem{lem}[thm]{Lemma}
\newtheorem{cor}[thm]{Corollary}

\newcommand{\pa}[1]{\left(#1\right)}
\newcommand{\cro}[1]{\left[#1\right]}

\newcommand{\E}{\operatorname{\mathbb{E}}}

\renewcommand{\P}{\operatorname{\mathbb{P}}}

\newcommand{\Lcal}{\mathcal{L}}
\newcommand{\Kcal}{\mathcal{K}}

\newcommand{\Pcal}{\mathcal{P}}

\newcommand{\Vcal}{\mathcal{V}}

\newcommand{\Zcal}{\mathcal{Z}}

\newcommand{\de}{\delta}
\newcommand{\si}{\sigma}
\newcommand{\ga}{\gamma}
\newcommand{\al}{\alpha}
\newcommand{\Va}{\textup{Var}}
\newcommand{\ovx}{\hat{x}}
\newcommand{\ovw}{\hat{w}}
\newcommand{\Var}{\textup{Var}}
\newcommand{\ptau}{p^{ \geq \tau}}
\newcommand{\ptaul}{p^{ < \tau}}

\newcommand{\ep}{\epsilon}

\newcommand{\Pb}{\mathbb{P}}

\newcommand{\ya}[1]{\textcolor{purple}{ya:[#1]}}

\title{Locally differentially private estimation of nonlinear functionals of discrete distributions}

\author{Cristina Butucea$^*$ \ 
and \, Yann Issartel\footnote{CREST, ENSAE, Institut Polytechnique de Paris, France. Yann.Issartel@telecom-paris.fr} }

\newcommand{\ovf}{\overline{F}}

\begin{document}

\maketitle

\begin{abstract}
We study the  problem of estimating non-linear functionals  of discrete distributions in the context of local differential privacy. The initial data $x_1,\ldots,x_n \in[K]$ are supposed i.i.d. and distributed according to an unknown discrete distribution $p = (p_1,\ldots,p_K)$. Only $\al$-locally differentially private (LDP) samples $z_1,...,z_n$ are publicly available, where the term 'local' means that each $z_i$ is produced using one individual attribute $x_i$. We exhibit privacy mechanisms (PM) that are sequentially interactive (i.e. they are allowed to use already published confidential data) or non-interactive.
We describe the behavior of the quadratic risk for estimating the power sum functional $F_\gamma= \sum_{k=1}^K p_k^{\ga}$, $\ga >0$ as a function of $K, \, n$ and $\alpha$. In the non-interactive case, we study two plug-in type estimators of $F_\ga$, for all $\ga >0$, that are similar to the MLE analyzed by Jiao \textit{et al.} \cite{Jiao_2017} in the multinomial model. However, due to the privacy constraint the rates we attain are slower and similar to those obtained in the Gaussian model by Collier \textit{et al.} \cite{collier2020estimation}. In the sequentially interactive case, we introduce for all $\ga >1$ a two-step procedure which attains the parametric rate $(n \alpha^2)^{-1/2}$ when $\ga \geq 2$.  We give lower bounds results over all $\alpha$-LDP mechanisms and all estimators using the private samples. 
\end{abstract}

\section{Introduction}

Information theoretic measures have become of utmost importance and extensively used in information theory, image processing, physics, genetics and more recently in machine learning and statistics. Such functionals of probability distributions are useful to design estimators, to choose most informative features in different algorithms, to test identity or closeness of distributions. Popular such measures are the power-sum, the entropy and more general R\'enyi entropies of discrete distributions.

%\ya{en tant que lecteur j'aurais du mal à comprendre les 2 paragraphes  ci-dessus qui sont peu concrets; pourquoi ne pas les remettre dans la section où ils appartiennent pour lier les choses similaires ?}

%\subsection{Power Sum Functional.}  
In this paper, we are interested in estimating the  power sum functional $F_\ga(p)$ of a discrete distribution $p=(p_1,...,p_K)$:
$$
F_\ga(p) = \sum_{k=1}^K p_k^\ga, \quad \text{ with power }\gamma \in (0,\infty).
$$

This instance of information measure has a tight connection with the famous Rényi entropy $H_{\ga}$  via the  formula $H_{\ga}= \frac{\log F_{\ga}}{1-\ga}$. 

In the statistical literature, smooth nonlinear functionals are often reduced via Taylor expansion to several functionals of the type $F_{\ga}$ for positive integer values of $\ga$ (see monographs like  e.g. \cite{GineNickl} ). For example, the entropy $H(p) = \sum_{j=1}^K p_j \log(1/p_j)$ of the probability distribution $p$ with finite support and probabilities bounded away from $0$ can be approximated via a Taylor expansion to a linear combination of $F_{\ga}$ for integer values of $\ga$ and estimated at parametric rate of $1/\sqrt{n}$.  %However, when the probabilities $p_j$ can be arbitrarily small, the derivative of $p_j \log(1/p_j)$ grows to infinity and this problem has been solved differently in this regime, in a series of papers as described later on. 

Symmetric functions of $p_1,\ldots,p_K$, i.e. functions of at most $K$ variables that are permutation invariant, play an important role in deep learning, e.g. \cite{Zaheer}, \cite{janossy}. Such functions can be written as polynomials of functionals $F_{\ga}$ for integer values of $\gamma$ and estimated using our procedures.

Another important application of such functionals is testing identity or closeness of distributions. For example, let us consider uniformity testing that is the null hypothesis is $H_0: p_j = 1/K$ for all $j$ from 1 to $K$ against the alternative hypothesis $p$ is not the uniform distribution. Suppose the distance measuring how far $p$ is from the uniform distribution is the Hellinger distance $\sum_{k=1}^K (\sqrt{p_k}-{1}/{\sqrt{K}})^2$. Therefore a test procedure will proceed by estimating this Hellinger distance using the sample and this involves estimating $F_{1/2}$ for $\ga=1/2$. Thus, an uniformity test procedure will be based on the estimator of this discrepancy, and similarly for identity or closeness tests. 
More generally, when the distance between two probability distributions is evaluated by a discrepancy or a distance, functionals $F_{\ga}$ naturally appear in their expression.
Note however that testing rates may differ from the estimation rates of the discrepancy as is the case for the $L_1$ distance where \cite{VV14} showed that several tests procedures must be aggregated in order to attain better rates for testing. We stress the fact that testing is a different problem from learning the functional.

\subsection{Plug-in Approach}  In the standard statistical setup (also called multinomial setting), the goal is to estimate the power sum functional $F_{\ga}$ based on $n$  i.i.d. samples $x_1,...,x_n$ following an unknown discrete distribution $p=(p_1,\ldots,p_K)$ with alphabet size $K$. A commonly-used approach to this problem is the plug-in approach, which  amounts to using an estimate $\hat{p}$ of the parameter $p$ in order to build an estimator $F_{\ga}(\hat{p})$ of the functional $F_{\ga}(p)$. The resulting plug-in estimator actually corresponds to the so-called maximum likelihood estimator (MLE) when the estimate $\hat{p}=(\hat{p}_k)_{k\in[K]}$ is defined as the empirical distribution $\hat{p}_k =\frac{1}{n} \sum_{i=1}^n\mathbbm{1}_{x_i = k}$ (i.e. the  average counts in the $k^{\textup{th}}$ box after binning).  This approach is not only intuitive and simple, but is also theoretically well grounded, as it is asymptotically efficient for finitely supported probabilities (finite $K$), and is non-asymptotically nearly optimal for possibly increasing $K$ \cite{Jiao_2017} (for details, see the related literature below). A natural question that we investigate in this paper is whether such plug-in type approach still performs well in a non-standard statistical setup where a constraint of privacy is imposed on the observed data. 

In the standard setup, estimation of power-sum functionals was  studied for the widely spread MLE estimator  in \cite{Jiao_2017} where the authors found that its maximal quadratic risk (also called worst case squared error risk) is  sub-optimal only by some  logarithmic factor. The estimation rates have then been  tightened to minimax optimal bounds using the best polynomial approximation of the power function \cite{Jiao_2015}, \cite{WuYang}. Following the chronology of this literature,  our study in a non-standard privacy setup starts with an analogue of this MLE, with a twofold purpose: (i) to highlight an important difference between the LDP setting and the standard (non private) setting: practical methods performing well in the non private case should not be systematically transferred to the private settings, a good estimator in the first setup  not being necessarily a good estimator in the latter setup;  (ii)  to show the regimes where  plug-in type estimators fail and the estimation problem is delicate, thus setting  benchmarks for future work on functional estimation in the LDP setting, as well as more involved  private settings.

\subsection{Differential Privacy}  

Keeping sensitive data $x_1,\ldots,x_n$ private  is a major concern in the modern area of Big Data. For example, $x_1,\ldots,x_n$ may be personal health or financial data of $n$ participants to a survey. \textit{Differential privacy} (DP) \cite{Dwork} has prevailed in the recent literature as a convenient approach to randomize samples with a control of the amount $\alpha >0$ of privacy introduced. %and of  the information contained in the original sensitive sample.
The randomized samples $z_1,\ldots,z_n$, also known as private samples, are provided to the statisticians who want to extract information on the underlying distribution $p$ of the initial data  $x_1,\ldots,x_n$.
\textit{Global or central DP} allows simultaneous treatment of the whole initial sample in order to produce the privatized random variables. In contrast,
\textit{Local differential privacy} (LDP)  is a stronger setup of privacy where no one has access to all sensitive data $x_1,\ldots,x_n$ (not even a trusted curator or third party  to handle the   privatization), but each individual $i$ has access to one $x_i$.

It is quite popular now that privacy, in particular local differential privacy, comes at the cost of slower rates of learning in many estimation problems. For example, estimation of the probability density is known to be achieved with slower rates in the minimax sense under $\alpha$-LDP constraints, see \cite{duchi2013}, \cite{WassermanZhou}, \cite{rohde2020geometrizing} and\\ \cite{butucea2020local}.  There is a rush nowadays to better understand when the loss is unavoidable and to show that such loss is optimal over all privacy mechanisms and procedures at hand.

In this LDP setup, there are major distinctions  according to the way the privacy mechanisms use the available information, that is non-interactive, sequentially interactive or fully interactive setups. When using {\it non-interactive} mechanisms, each private sample $z_i$ is produced by a privacy mechanism $Q_i$ that has access only to one sensitive sample $x_i$. They are arguably the simplest mechanisms. A richer class of mechanisms  are the so-called \textit{sequentially interactive} mechanisms where each agent $i$ is allowed to incorporate into its private mechanism $Q_i$ the already privatized data $z_1,\ldots,z_{i-1}$ of other agents, along with $x_i$. It is known e.g.  \cite{butucea2020interactive}  that sequentially interactive methods  attain much  faster rates than the non-interactive mechanisms in some inference problems. This has also been proved for identity testing of discrete distributions in \cite{BerrettButucea}. Finally, the {\it fully interactive} privacy mechanisms are allowed to use one $x_i$ several times, together with all publicly available randomized $z$'s. In this case, the natural question is how many randomizations are necessary in order to acquire the desired amount of information. It has been proved that there exist separations  between fully and sequentially interactive procedures, e.g. \cite{josepha} and \cite{josephb}, as well as separations between fully interactive mechanisms and global DP, e.g. \cite{Chan}. We do not consider fully interactive mechanisms in the current paper (neither in the upper, nor in the lower bounds).

Understanding the relative power of interactivity is a crucial question in the LDP setting. In this line of work, we study the estimation of the power sum functional $F_{\ga}$ under the constraints that the available data $z_1,\ldots,z_n$ stem from non-interactive and sequentially interactive privacy mechanisms, respectively.

\subsection{Contributions} 

In the $\alpha$-LDP setting % where we only have access to privatized versions $z_1,\ldots,z_n$ of the original data $x_1,\ldots,x_n \overset{\text{i.i.d.}}{\sim} p$, 
we propose three estimators of the power sum function $F_{\ga}(p)$ based on the $n$ randomized observations. The performance of any estimator $\hat{F}_{\ga}$ is controlled by proving upper bounds on its non-asymptotic quadratic risk $\E[(\hat{F}_{\ga} - F_{\ga})^2]$, with an explicit dependence on the parameters of the problem: the power $\ga$, the alphabet size $K$, the sample size $n$ and the amount of privacy $\alpha$. 

Our first contribution is a tight characterization of the quadratic risk of the \textit{plug-in estimator}, an analogue of the plug-in MLE (discussed above). Unfortunately, its risk grows rapidly with the alphabet size $K$, thus showing that the plug-in approach is by far not optimal for large $K$.  This contrasts with the good performance of the plug-in in the standard statistical setup \cite{Jiao_2017}. Thus, good estimators in the standard setup are not necessarily to be used as such in the LDP setting. 

Our second contribution is a correction of this plug-in estimator by truncating the small probabilities $p_k$. The induced procedure, called \textit{thresholded estimator} in the sequel, performs significantly better than the plug-in estimator, typically when the alphabet size $K$ is large.
We emphasize that this  improvement is important, since the risk of this thresholded estimator is (almost) independent of the support size $K$ for large $K$, unlike the risk of the plug-in estimator. It is therefore different from the literature on functional estimation (in the standard setting)  where the improvements in the risk of plug-in estimators are often of magnitude of logarithmic factors. 

The privatized data $z_1,\ldots,z_n$ used by both our plug-in and thresholded estimators are generated by a simple Laplace non-interactive mechanism. In contrast, our third contribution is a two-step procedure based on a sequentially  interactive mechanism. The definition of this two-step procedure heavily relies on the plug-in estimator, and thus can be seen as a refinement of the plug-in approach. Such a sequentially interactive method was studied in \cite{butucea2020interactive} for the particular power sum functional $F_{2}$ ($\ga =2$), also called quadratic functional, in a continuous setup  (where the probability distribution is a smooth function on $[0,1]$). By allowing to encode information from previous observations $z_1,\ldots,z_{i-1}$ into new released data $z_i$, % knowledge of the probabilities using  the non-interactive mechanism for $z_1,\ldots,z_{[n/2]}$ in the new released data $z_{[n/2]+1},..., z_n$,
this sequentially interactive procedure  achieves faster rates than the thresholded estimator, when $K$ is large and $\ga >1$, though this improvement is only of a logarithmic factor. Unfortunately, this sequentially interactive procedure is only defined for $\ga >1$, and has slower rates than the thresholded and plug-in estimators when $K$ is small.  Accordingly, none of our three estimators is overall better than the others. Table~\ref{Rates} presents upper bounds on the maximal quadratic risks of these three estimators. The choice of estimators therefore depends on regimes (i.e. values of $K, \ga$), and our fastest rates attained by some combination of these three estimators  are written in Corollary \ref{mainTHM}. 

We finally give lower bounds on the maximal quadratic risk, over all estimators and all (non-interactive and sequentially interactive) privacy mechanisms. These lower bounds are optimal for $\ga \geq 2$ as they match our fastest rates. Unfortunately for $\ga \in (0,2)$, a gap of a factor $K^{\ga}$ remains between our lower bounds and fastest rates.   All proofs are in \cite{butIss2021nips}. 

\begin{table}
  \caption{ Upper bounds on maximal quadratic risks of three procedures}
  \label{Rates}
  \centering
  \begin{tabular}{lccc}
    \toprule
    & $\gamma \in (0,1)$ &  & $\ga > 1$ \\
    \midrule
   % & & \\
  %  \textit{Non-Interactive PM:} & & \\
  %  & & &\\
  $ \begin{array}{l} \text{Plug-in estimator} \\
   (\textit{Non-Interactive PM})  \end{array}$
    &  $\frac{K^2}{(\al^2 n)^\ga} $
    &
    &  $\frac{K^2}{(\al^2 n)^\ga} 
    + \frac{ K^{3 - 2\ga} \vee 1 }{\al^2 n} \ \, \, := R_1$  \\
    & & \\
  $ \begin{array}{l} \text{Thresholded  estimator } \\
   (\textit{Non-Interactive PM}) \end{array}$
    &  $  K ^{2 (1- \ga)} \wedge \frac{K^2 }{(\al^2  n)^{\ga }} $
    &
    &  $\min \{ R_1,R_2 \} . \log(Kn)^{\ga} $ \\
   % \begin{array}{l}   \min \left\{ \frac{\left(\log(Kn)\right)^{\ga-1}}{(\al^2 n)^{\ga -1}} + \frac{1 }{\al^2  n } ,\right. \\
   % \\
   %  \left.  \frac{K^2 \left(\log(Kn)\right)^{\ga}}{(\al^2 n)^{\ga }}  + \frac{ K^{3-2\ga} \vee 1 }{\al^2  n } \right\}
   % \end{array} $\\
    & & \\
    \midrule
    %& & \\
  % \textit{Interactive PM:} & & \\
  % & & &\\
   Two-step procedure (\textit{Interactive PM}) &   
   &
   &  $ \frac{ 1 }{(\al^2 n)^{\ga - 1}}
    + \frac 1{\al^2 n} \ \, \,  := R_2$  \\
    \bottomrule
  \end{tabular}
\end{table}

%%%%%%%%%%%
% new section
%%%%%%%%%%%%%

\subsection{Related Literature }
 
 To the best of our knowledge, estimating non linear functionals under global DP has been considered by \cite{Acharya2018}. The authors estimate the entropy, the support size and the support coverage of a discrete distribution in the context of global differential privacy. In each setup they show that the cost of privacy is relatively small when compared to the standard (non private) setup. Their upper bounds allow to quantify the amount of privacy that we can have without deteriorating the estimation rates.

We study for the first time  the estimation of the nonlinear power-sum functional $F_\ga$ for any real $\ga >0$ in the context of LDP.  The case $\ga=2$ has been considered in the LDP setup,  for smooth distributions in \cite{butucea2020interactive}, and  for testing discrete distributions in \cite{BerrettButucea}.  In \cite{butucea2020interactive}, the association of a plug-in estimator and Laplace mechanism is optimal among all non-interactive mechanisms, whereas a sequentially interactive procedure improves dramatically the minimax rates. We will also try to further understand when such phenomena hold for different values of $\ga$, in particular when the functional $F_{\ga}$ is less smooth for $\ga <2$.

In non-private settings, it has been shown different behaviors for estimating $F_\gamma$ according to
the observations scheme, namely faster rates are attained in the multinomial setup \cite{Jiao_2017} than in the Gaussian vector  model \cite{collier2020estimation} $-$ see  details below. We will see that due to the LDP setup our rates are similar to those proved in the Gaussian vector  model by \cite{collier2020estimation}, though our plug-in type estimator is an analogue of the MLE analyzed in the multinomial  setting by \cite{Jiao_2017}. 

\textit{Standard (or multinomial)  setting: } When the r.v. $x_1,...,x_n$ are observed,
\cite{Jiao_2017} show that the maximal quadratic risk of the maximum likelihood estimator (MLE) in estimating $F_{\ga}$, $\ga >0$, is
\[%\begin{equation*}
   \frac{K^2}{n^{2 \ga}} \mathbf{1}_{ \ga \in (0,1)} + \frac{K^{2(1-\ga)}  }n \mathbf{1}_{\frac 12 < \ga < 1 }
+ n^{-2(\ga-1)} \mathbf{1}_{ 1<\ga < \frac{3}{2} } + \frac 1n \mathbf{1}_{ \ga \geq \frac{3}{2}}.
%\end{equation*}
\]
In this model, the MLE of $p_k$ is the average of the $\mathbf{1}_{x_i=k}$ over $i$ from 1 to $n$. These results entail that the MLE achieves the minimax rate $n^{-1}$ when $\ga \geq 3/2$. However, when the regularity $\ga$ decreases below $3/2$, the rates of the MLE  get slower. In particular, the MLE does not achieve the minimax rates of estimation of $F_{\ga}$ when $0 <\ga< 3/2$. For $0 <\ga<1$, the difficulty of estimation increases, especially if the alphabet size $K$ is large, which may even prevent from a consistent estimation of $F_{\ga}$ (i.e. from the risk of the MLE to converge to zero). More precisely, the MLE consistently estimates $F_{\ga}$, $0<\ga <1$, if and only if, the sample size is large enough to satisfy $n \gg K^{1/\ga}$.

For the values of $\ga < 3/2$,  the minimax rates in estimating $F_{\ga}$ roughly correspond to the performance of the MLE (written above) with $n$ replaced by $n\log n$ \cite{Jiao_2015}. These faster rates have been attained through a different estimator that uses best polynomial approximation of the functional $p_k^\ga$ for small values of $p_k$, which allows to gain a log factor in the bias. This method has been largely employed in the literature since the pioneer work by \cite{Lepski}, in the Gaussian white noise model. The best polynomial approximation technique has been used for estimating the $\mathbb{L}_1$ norm in a Gaussian setting \cite{CaiLow}, the $\mathbb{L}_r$ norm in the Gaussian white noise model \cite{Han_a}, the entropy of discrete distributions \cite{WuYang}, the entropy of Lipschitz continuous densities \cite{Han_b}.
In the multinomial setting, \cite{Acharya} studied the estimation of the R\'enyi entropy of order $a>0$ and \cite{Fukuchi} extended these results  to more general additive functionals of the form $\sum_{k=1}^K \varphi(p_k)$, for 4-times differentiable functions $\varphi$.

%Best polynomial approximation (e.g. see [11] to [14] \ya{Réf en dur}) achieves  the minimax rate in the non-private case. For the power sum, we would expect rates faster by a logarithmic factor than ours. (a detailed description of the method is given in the rebuttal to the last reviewer). This method will be the object of further investigation mainly because the gain of a logarithmic factor did not seem to us the most pressing question in the LDP estimation of $F_{\ga}$, as there are still important missing pieces in the comprehension of the problem, mainly proving lower bounds that do not suffer from a factor $K^{\ga}$. 

\textit{Gaussian vector setting:} the observations are $y_i = \theta_i + \ep \cdot  \xi_i$, $i=1,\ldots,K$ where $(\theta_1,\ldots,\theta_K)$ are unknown parameters, $\ep>0$ is known, and $\xi_i$ are i.i.d. standard Gaussian random variables. The authors of \cite{collier2020estimation} show that there exists an estimator based on best polynomial approximation with the following error bound $-$ see (14) in the proof of their Theorem 1 $-$
\begin{equation}\label{GaussianSetting:1}
  \E \cro{ \pa{\hat{F}_{\ga} - F_{\ga}}^2}  \lesssim \frac{\ep^{2\ga}K^2}{\log^{\ga}(K)} + \frac{\ep^2}{\log(K)} \left( \sum_{i=1}^K \theta_i^{\ga -1} \right)^2 \mathbbm{1}_{\ga >1}
 + \ep^2 \sum_{i=1}^K \theta_i^{2\ga -2}\mathbbm{1}_{\ga >1}
 \end{equation}for any real $\ga >0$.
When $\ga \in(0,1)$,  this bound is equal to $\ep^{2\ga}K^2 / \log^{\ga}(K)$ and turns out to be minimax optimal \cite{collier2020estimation}.
 Besides, the authors show that better rates can be achieved if $\ga$ is an integer. Namely, there exists an estimator with the following error bound, see Theorem 2 in \cite{collier2020estimation},
\begin{equation}\label{GaussianSetting:2}
  \E \cro{ \pa{\hat{F}_{\ga} - F_{\ga}}^2}  \lesssim \ep^{2\ga}K + \ep^2 \| \theta \|_{2\ga -2}^{2\ga -2}
 \end{equation}for any integer $\ga\geq 1$. In this special case of an integer $\ga$, the rate \eqref{GaussianSetting:2} is achieved by a simple estimator that has no bias.

%%%%%%%%%%%%%%%%%%%%%%%%%%%%%%

\section{Estimators and results}  

We aim at estimating the power sum functional $F_{\ga}(p) =  \sum_{k=1}^K p_k^{\ga} $, for $\ga \in (0,\infty)$, in the LDP setting where we only have access to privatized versions $z_1,\ldots,z_n$ of the sensitive original data $x_1,\ldots,x_n \overset{\text{i.i.d.}}{\sim} p$. 
The sensitive random variables $x_1,...,x_n$  are i.i.d. distributed according to the discrete density model $ \Pb ( x_i = k ) = p_k$  for $k= 1, \ldots, K$
where the unknown parameter $p = (p_1,\ldots,p_K)$ belongs to  the set  $\Pcal_K = \{ p \in [0,1]^K \, : \sum_{k=1}^K  p_k =1\}$ of discrete distributions with alphabet size $K$. 
The $x_1,\ldots,x_n$ are not publicly available; instead they are used as inputs into a privacy mechanism (PM), in order to produce available sanitized observations $z_1,...,z_n$.
We use the sanitized observations to estimate the functional $ F_\ga$. 

Notation: The symbol $x \lesssim_\ga y$ means that the inequality $x \leq C_{\ga} y$ holds for some constant $C_{\ga}$ depending only on $\ga$.
We denote  $\min \{ x,y\}$ by $x \wedge y$, and  $\max \{ x,y\}$  by $x \vee y$.

%In this section, we first introduce a non-interactive PM, then we present two estimators of $F_{\ga}$ that use the private sample produced by this  PM. The first estimator is a plug-in estimator which, unfortunately, has a risk that grows rapidly with the alphabet size $K$. So, with the objective of reducing the effect of $K$ on the risk, we introduce  a refinement of the plug-in estimator, called thresholded estimator.

%\ya{truc}We first consider two non-interactive methods, where we privatize the initial sample via a standard Laplace mechanism. The first estimator is an analogue of the widely-used MLE \cite{Jiao_2017} discussed in the introduction. It uses an unbiased estimation of the vector $p$, by averaging the privatized samples in each bin $k=1,\ldots,K$. The resulting estimator of $F_\ga$ is actually not a MLE in this privatized context, hence we call it a plug-in estimator. The convergence rates of this plug-in are slower than those of the MLE \cite{Jiao_2017}, although both estimators have similar definitions. In fact, the privatized Laplace setup is close to the Gaussian vector model \cite{collier2020estimation}, and our plug-in estimator achieves similar rates to those in \cite{collier2020estimation}.
 
\subsection{LDP setup}
 
The privacy mechanism PM also known as channel or randomization, is submitted to the constraint that it is $\alpha$-locally differentially private (LDP) for some $\alpha >0$. This means that the PM generates the private samples $z_i$  using a conditional distribution $Q_i(\cdot | x_i, z_1,...,z_{i-1})$ such that 
\begin{equation}\label{def:alphaLDP}
 \sup_{z_1,...,z_{i-1}, x_i, x_i'} \frac{Q_i(\cdot | x_i, z_1,...,z_{i-1})}{Q_i( \cdot | x'_i, z_1,...,z_{i-1})} \leq e^\alpha, \quad \mbox{for all } i=1,...,n,
\end{equation}
with the convention that  $\{z_1,\ldots,z_{i-1}\}$ is the empty set for $i=1$. 
These PM are called sequentially interactive, as each $Q_i$ is allowed to use  previously published samples $z_1,...,z_{i-1}$. We are also interested in the sub-class of non-interactive  PM which are not allowed to use previous released data. A non-interactive PM generates each $z_i$ only accessing $x_i$,  via a conditional distribution of the form $Q_i(\cdot | x_i )$. We assume from now on that the level $\alpha$ of privacy belongs to $(0,1)$ and satisfies $\al^2n \geq 1$.

\subsection{Non-interactive privacy mechanism (NI PM)} 
We introduce a non-interactive PM, denoted by $ Q^{(NI)}$. Given the original data $x_i,$ individual $i$ generates a random vector $z_i =(z_{i1},\ldots,z_{iK})$ using the Laplace non-interactive privacy mechanism $Q^{(NI)}$  defined by
\begin{equation}\label{NIprivacyMechanism}
 Q^{(NI)}: \qquad   z_{ik} = \mathbf{1}_{\{x_i = k\}} +\frac{\si}{\alpha} \cdot w_{ik}, \qquad k=1, \ldots,K,
\end{equation}
where the $w_{ik}$ are i.i.d. Laplace distributed with density $f^w(x)= \frac{1}{2} \exp(-|x|).$ Note that $w_{ik}$ are all centered, with variance $2$. Setting $\si=2$, one can readily check that the channel $Q^{(NI)}$  above is an  $\alpha$-LDP non-interactive mechanism, see \cite{BerrettButucea} and the references therein. 
We denote the mean of the privatized observations in the $k^{\textup{th}}$-box by $\hat z_k = \frac 1n \sum_{i=1}^n z_{ik}$.

\subsection{Plug-in estimator based on NI PM} 
The first estimator we introduce is an analogue of the  MLE \cite{Jiao_2017} discussed in the introduction. It uses an unbiased estimation of the parameter $p$, by averaging the privatized samples in each bin, thus using the $\hat{z}_k$,  $k=1,\ldots,K$.
The resulting estimator $\hat{F}_{\ga}$ of $F_\ga$ is actually not a MLE in our privatized context, hence we call it plug-in estimator. The plug-in  estimator uses the privatized data produced by the NI PM above, and then estimates separately each term  $F_{\ga }(k) = p_k^\gamma$ of the functional $F_{\ga }= {\sum_{k=1}^K} F_{\ga }(k)$  as follows:   
\begin{equation}\label{nonInteractive:procedure}
\hat{F}_{\ga} = \sum_{k=1}^K \hat F_\ga (k) \enspace, \mbox{ with }
 \hat F_\ga (k) =  \left(T_{[0,2]}\cro{\hat z_{k}} \right)^\ga \enspace,
\end{equation}
where $T_{[0,2]}\cro{\cdot}$ is the clipping operation defined by $T_{[0,2]}\cro{y} =  \pa{y \vee 0} \wedge 2 $. As $\hat{z}_k$ is an unbiased estimator of $p_k$ with fluctuations of order 
$(\sqrt{ p_k(1-p_k)/n}) +  (\sigma/\sqrt{ \al^2 n })$, the quantity $\tau = c / \sqrt{\al^2  n}$, where $c\geq 1$ is a numerical constant,  can be seen as a threshold that is just above the noise level in the available data.
%We separate the coordinates of $p$ accordingly, denoting by $\Kcal_{> \tau} = \{ i\in [K]: \, p_i > \tau \}$ (respectively $\Kcal_{\leq \tau} = \{ i\in [K]: \, p_i \leq  \tau \}$) the set of coordinates above $\tau$ (resp. below $\tau$).
Write $\ptau$  (respectively $\ptaul$) the vector containing the thresholded values of $p$: $p_k \cdot \mathbf{1}(p_k \geq \tau)$ (respectively $p_k \cdot \mathbf{1}(p_k < \tau)$) for  $k\in [K]$.

\begin{thm}\label{thm:bias:nonInterac:summary:new}
For any $\ga >0$ and  $p\in \Pcal_K$, the quadratic risk of the estimator (\ref{nonInteractive:procedure}) is bounded by
\begin{equation}\label{risk:tauPetit:thm}
  \E_p \left[( \hat F_\gamma - F_\gamma)^2 \right] \lesssim_\ga \frac{K^2}{(\al^2  n)^{\ga}}
    + \mathbf{1}_{\{  \ga \geq 2\}}\frac{\| \ptau \|^{2(\ga - 2) }_{\ga -2}}{(\al^2  n )^2 }
    + \mathbf{1}_{\{  \ga \geq 1\}}\frac{\| \ptau \|^{2\ga -2}_{2\ga -2}}{\al^2  n } \enspace,
\end{equation}
and thus is uniformly bounded over all $p\in \Pcal_K$ by
\begin{equation*}
    \sup_{p\in \Pcal_K} \E \cro{ (\hat{F}_{\ga} - F_{\ga})^2 } \lesssim_{\ga}  \frac{K^2}{(\al^2  n)^{\ga }}
%      + \mathbf{1}_{\{  \ga \geq 2\}}\frac{1 \vee K^{2(3-\ga)}}{((\al^2 \wedge 1) n )^2 }
      + \mathbf{1}_{\{  \ga \geq 1\}}\frac{( \tau^{-1} \wedge K)^{3-2\ga} \vee 1 }{\al^2 n } \enspace.
\end{equation*}
\end{thm}

This upper bound on the quadratic risk of the plug-in estimator grows quadratically with the alphabet size $K$, even for arbitrarily  large $K$. This poor performance is unfortunately inherent to  the plug-in estimator, as shown by the lower bound  in Proposition \ref{borne:inf:plugin:petit:pk}.  Accordingly, we will correct the plug-in estimator in the next section to reduce the effect of the alphabet size $K$ on the risk.

\begin{prp}\label{borne:inf:plugin:petit:pk}
For any $\ga >0$, $\ga \neq 1$,  and integer $K \geq 2$, the maximal quadratic risk of $\hat{F}_{\ga}$ is asymptotically  bounded from below by
\begin{equation*}
   \sup_{p\in \Pcal_K}  \E \cro{ (\hat{F}_{\ga} - F_{\ga})^2 } \gtrsim_{\ga}  \frac{K^2}{(\al^2  n)^{\ga}}  \qquad \textup{as} \ \, n\rightarrow \infty \enspace.
\end{equation*}
\end{prp} 

The maximal quadratic risk of the plug-in estimator has therefore a quadratic dependence on the alphabet size $K$, and thus it is necessary and sufficient to have $n \gg K^{2/\ga} \al^{-2}$ initial observations for this estimator to be consistent (in the sense that its risk converges to zero). This rate is slower than that of the MLE  in the standard setup without privacy constraint \cite{Jiao_2017}, though both our plug-in estimator and the MLE have  similar definitions based on the plug-in principle. Indeed,  the MLE of \cite{Jiao_2017} is less sensitive to the support size $K$ than our plug-in estimator, as the MLE is consistent regardless of the alphabet size $K$, as soon as $\ga >1$. In addition, when $\ga\in(0,1)$, the MLE is consistent if, and only if $n \gg K^{1/\ga}$. This better performance of the MLE can be explained by the fact that \cite{Jiao_2017} benefits from direct observations $\hat{x}_k := \frac{1}{n} \sum_{i=1}^n \mathbf{1}_{\{x_i = k\}}$ that have non-homogeneous fluctuations, i.e. a variance $\textup{Var}(\hat{x}_k)= p_k(1-p_k)/n$ that scales with  the signal $p_k$. By contrast, our plug-in estimator employs privatized observations  $\hat{z}_k$ resulting from the Laplace PM in \eqref{NIprivacyMechanism}, which have nearly homogeneous fluctuations with a variance  $\textup{Var}(\hat{z}_k)$ that scales  with $\tau^2 \asymp (\al^2n )^{-1}$ regardless of $k\in[K]$.  Thus, both situations are different for small $p_k$, which typically occurs when $K$ is large since $\sum_{k=1}^K p_k=1$.

\textit{Link with the Gaussian vector model.} The observations released by the Laplace PM in \eqref{NIprivacyMechanism} are in fact close to the observations of the Gaussian vector model \cite{collier2020estimation} and, as a consequence, our plug-in estimator achieves similar rates to those obtained in \cite{collier2020estimation}. Indeed, our observations $(\hat{z}_k)_{k\in [K]}$  can be seen as homoskedastic random variables having common characteristics to the Gaussian random variables in \cite{collier2020estimation} with variance  $\epsilon^2_0= 1/(\al^2 n)$. Let us check that our upper bound \eqref{risk:tauPetit:thm} implies, up to log factors, the same bound as \eqref{GaussianSetting:1} from \cite{collier2020estimation}. (The logarithmic gap comes from the fact  that we do not use the best polynomial approximation, unlike \cite{collier2020estimation} $-$ see the introduction for details.)    The square root of the second term $\| \ptau \|^{\ga -2}_{\ga -2}/(\al^2  n)$ in \eqref{risk:tauPetit:thm} is bounded from above by
\begin{equation*}
     \frac{ \sum_{k=1}^K p_k^{\ga -2} \mathbbm{1}_{\{p_k > \tau\}}}{\al^2  n }
     \leq \sum_{k=1}^K \frac{  p_k^{\ga -1} }{\sqrt{\al^2  n} }  \frac{ p_k^{-1} \mathbbm{1}_{\{p_k > \tau\}} }{\sqrt{\al^2 n} }
     \leq \sum_{k=1}^K \frac{  p_k^{\ga -1}  c^{-1}}{\sqrt{\al^2  n} }
     =  \frac{  \| p \|^{\ga -1}_{\ga -1} c^{-1} }{\sqrt{\al^2  n} } 
\end{equation*}which corresponds to the square root of the second term of (\ref{GaussianSetting:1}), up to a log factor. Thus, our bounds are similar to those in \cite{collier2020estimation} for a specific variance $\epsilon^2_0$. Note the difference that we impose the constraint $p\in \Pcal_K$, whereas  \cite{collier2020estimation}  consider all $p\in \mathbb{R}^K$.

%%%%%%%%%%%%%%%%%%%%
% new section
%%%%%%%%%%%%%%%%%%%

\subsection{Thresholded plug-in estimator based on NI PM}
\label{thresholdsubsection:thm}

Our objective is to reduce the effect of the alphabet size $K$ on the risk of the plug-in estimator. A natural solution  is  to not estimate the small probabilities $p_k$, as their contribution to the functional $F_{\ga}$ is relatively weak, while they may add much fluctuations to the plug-in estimator $\hat{F}_{\ga}= \sum_{k=1}^K \hat{F}_{\ga}$.  
Accordingly, we introduce a refinement of the plug-in estimator, called \textit{ thresholded plug-in estimator} $\overline{F}_{\ga}$, which does not estimate all components $p_k^{\ga}$ of  the sum $F_{\ga}=\sum_k p_k^{\ga}$. Similar to the plug-in estimator, this estimator $\overline{F}_{\ga}$ uses the private sample $(z_{ik})_{i \in[n], k\in[K]}$ produced by the non-interactive PM in  \eqref{NIprivacyMechanism}.

The  definition of $\overline{F}_{\ga}$ is split in two cases, $\ga \in(0,1)$ and $\ga >1$. When $\ga \in(0,1)$, we simply set $\overline{F}_{\ga} := \mathbbm{1}_{K\leq \tau^{-1}} \hat{F}_{\ga}$, meaning that $\overline{F}_{\ga}$ is equal the trivial estimator $0$ if $K\geq \tau^{-1}$, and  to the plug-in estimator $\hat{F}_{\ga}$  otherwise. %We thus say that for large $K\geq \tau^{-1}$ there is no point in estimating the functional $F_{\ga}$, while for small $K \leq \tau^{-1}$  the plug-in estimator is fine.  

For $\ga > 1$, the thresholded estimator $\overline{F}_{\ga}$ uses the plug-in estimator only on the significant components $p_k^{\ga}$ of the sum $F_{\ga}=\sum_k p_k^{\ga}$, via a truncation of the small probabilities $p_k$. This two-step procedure first detects the significant probabilities $p_k$ that are above the threshold $\tau =c /\sqrt{\al^2 n}$, and then estimates the part of the functional $F_{\ga}$ induced by those $p_k$. Assume that the sample size is $2n$ for convenience, and split the data in two samples $x^{(1)} =(x_1^{(1)},\ldots,x_n^{(1)})$ and $x^{(2)} =(x_1^{(2)},\ldots,x_n^{(2)})$. The individuals owning the data $z^{(s)}$, $s=1,2$, use the non-interactive mechanism \eqref{NIprivacyMechanism} 
which generates $z_i^{(s)}=(z_{i1}^{(s)},\ldots,z_{iK}^{(s)})$ for $i=1,\ldots,n$. Denote the two sample means of the $k^{\textup{th}}$ bin by  $\hat z_{k}^{(s)} = \frac{1}{n}\sum_{i=1}^n z_{ik}^{(s)}$, $s=1,2$.
We use the $(\hat{z}_{k}^{(1)})_k$ to detect large values of the underlying probabilities as follows. For each $k\in[K]$, if $\hat z_{k}^{(1)}$ is strictly smaller than the the empirical  threshold $\hat{\tau}$,
\begin{equation*}
  \hat{\tau} :=  192 \si \sqrt{ \frac {\log(Kn) }{ \al^2n} }  \enspace,
\end{equation*} then we do not estimate  $F_{\ga}(k)$. The threshold $\hat{\tau}$ is chosen larger than $\tau$ by a logarithmic factor, only to have high probability concentration for the $(\hat{z}_k^{(1)})_k$ around their means $(p_k)_k$. Otherwise, when $\hat z_{k}^{(1)} \geq \hat{\tau}$, we estimate $F_{\ga}(k)$ using the same plug-in estimator as in \eqref{nonInteractive:procedure}. This gives the following estimator $\overline{F}_{\ga}$ of $F_{\ga}$, $\ga >1$,
\begin{equation}\label{estim1:threshold}
   \overline{F}_{\ga} = \sum_{k=1}^K \overline F_\ga (k) \enspace, \mbox{ with }
 \overline{F}_{\gamma}(k) := \left(T_{[0,2]}\cro{\hat z_{k}^{(2)} \cdot \mathbf{1}_{\hat z_{k}^{(1)} \geq \hat{\tau}}} \right)^{\ga} \enspace,
\end{equation}where the $(\hat{z}_{k}^{(2)})_k$ are used to estimate the functional. The independence between $(\hat{z}_{k}^{(1)})_k$ and $(\hat{z}_{k}^{(2)})_k$ allows us to avoid cumbersome statistical dependencies between the two stages of the procedure.

\begin{thm}\label{threshold:thm} For any integers $K,n$ satisfying $ n \geq 2\log(K)$, the quadratic risk of $\overline{F}_{\gamma}$ is uniformly bounded over all $p\in \Pcal_K$ by \\
\textup{1)} $\ga\in(0,1)$:
\begin{equation*}
     \sup_{p\in \Pcal_K}  \E \cro{ (\overline{F}_{\ga} - F_{\ga})^2 }  \lesssim_{\ga}  K ^{2 (1- \ga)} \wedge \frac{K^2}{(\al^2 n)^{\ga }} \enspace.
\end{equation*}
\textup{2)} $\ga >1$:
\begin{align*}
    \sup_{p\in \Pcal_K}  \E \cro{ (\overline{F}_{\ga} - F_{\ga})^2 } &\lesssim_{\ga}   \bigg{(} \frac{\left(\log(Kn)\right)^{\ga-1}}{(\al^2 n)^{\ga -1}} + \frac{1 }{\al^2  n }\bigg{)} \wedge \bigg{(} \frac{K^2 \left(\log(Kn)\right)^{\ga}}{(\al^2 n)^{\ga }}  + \frac{ K^{3-2\ga} \vee 1 }{\al^2  n } \bigg{)}.
\end{align*} 
\end{thm}
Compared to the plug-in estimator performance,  the quadratic risk of the thresholded estimator $\overline{F}_{\ga}$ is much less sensitive to the alphabet size $K$ when $\ga >1$ and $K$ is large. Indeed, our bound on this risk is a minimum between two error terms where the left-term only depends on $K$ logarithmically. However, the new logarithmic factors that come from our use of high probability concentration  are not satisfactory.  We remove them in the next section, replacing the current non-interactive PM with a  sequentially  interactive PM.

%new section

\subsection{Two-step procedure based on sequentially interactive PM}

A sequentially  interactive PM is allowed to use the prior released (sanitized) data to encode our present knowledge  in the new released data. The idea is to rewrite the functional $F_{\ga}$ as $ \sum_{k=1}^K p_k \cdot F_{\ga - 1}(k)$, so that we  first compute the plug-in estimator $\hat{F}_{\ga-1}(k)$ of $F_{\ga-1}(k)$, and then build on $\hat{F}_{\ga-1}$ to estimate the functional $F_{\ga}$. In this two step procedure,  half of the sample is released via the Laplace mechanism \eqref{NIprivacyMechanism} and  is used for computing  $\hat{F}_{\ga-1}(k)$,  then the other half is released through a sequentially interactive PM encoding the information from $\hat{F}_{\ga-1}(k)$.  This procedure is only applicable for $\ga >1$ since it requires to compute $\hat{F}_{\ga-1}$. The sequentially interactive PM we consider here is similar to the ones studied for the specific functional $F_{2}$ ($\ga =2$) in  continuous setting \cite{butucea2020interactive} and identity testing \cite{BerrettButucea}.

%We consider the following sequentially interactive channel $Q^{(I)}$.
Assuming that the sample size is $2n$ for convenience, we split the data in two groups $x^{(1)} =(x_1^{(1)},\ldots,x_n^{(1)})$ and $x^{(2)} =(x_1^{(2)},\ldots,x_n^{(2)})$. The individuals owning the data $x^{(1)}$ use the non-interactive mechanism (\ref{NIprivacyMechanism}), i.e.
\begin{equation*}%\label{NIprivacyMechanism:inter}
 Q^{(NI)}: \qquad   z_{ik}^{(1)} = \mathbf{1}_{\{x_i^{(1)} = k\}} +\frac{\si}{\alpha} \cdot w_{ik}, \qquad k=1, \ldots,K,
\end{equation*}
which generates $z_i^{(1)}=(z_{i1}^{(1)},\ldots,z_{iK}^{(1)})$ for $i=1,\ldots,n$. Denote this first sample by $z^{(1)}=(z_1^{(1)},\ldots,z_n^{(1)})$. These sanitized data allow us to estimate $F_{\ga -1}(k) = (p_k)^{\gamma-1}$ using the plug-in estimator  (\ref{nonInteractive:procedure}), i.e.
\begin{equation*}%\label{estim1}
    \hat{F}_{\gamma-1}^{(1)}(k) = \left(T_{[0,2]}\cro{\frac{1}{n}\sum_{i=1}^n z_{ik}^{(1)}} \right)^{\ga -1} =: \left(T_{[0,2]}\cro{\hat z_{k}^{(1)}} \right)^{\ga-1} \enspace.
\end{equation*}
We then design an estimator of $\sum_{k=1}^K p_k \hat{F}_{\gamma-1}^{(1)}(k)$ which can be seen as a proxy of $F_{\ga}$. This is possible with the following sequentially interactive mechanism that encodes prior information $\hat{F}_{\gamma-1}^{(1)}$ in the released data $z^{(2)}_{i}$:
\begin{equation*}%\label{InteracMec:def}
 Q^{(I)}: \qquad  z^{(2)}_{i} =
         \pm z_\alpha, \mbox{with probability  } \frac{1}{2} \left(1 \pm \frac{ 1  }{ z_\alpha } \hat{F}_{\gamma-1}^{(1)}(x_i^{(2)}) \right)
\end{equation*}
where $z_\alpha = {2^{\gamma-1}} \frac{e^{\al}+1}{e^{\al}-1}$, and with the following sequentially interactive estimator 
\begin{equation*}%\label{def:estima:inter}
\widetilde{F}_{\ga} = \frac{1}{n} \sum_{i=1}^n z_i^{(2)} \enspace.
\end{equation*}
Denoting the second sample by  $z^{(2)}=(z_1^{(2)},\ldots,z_n^{(2)})$, it is easy to see that  the privatized sample $(z^{(1)}, \, z^{(2)})$  satisfies \eqref{def:alphaLDP} and thus is $\alpha-$LDP .  % This privatized sample is interactive in the sense that the second half also uses the first half of private samples in the process of randomizing locally each individual data $x_i$. 

\begin{thm}\label{thm:bias:summary:interac}
For any  $\ga >1$ and $p\in \Pcal_K$, the quadratic risk of the sequentially interactive estimator $\widetilde{F}_{\ga}$ is bounded by
\begin{equation*}
    \E \cro{ (\widetilde{F}_{\ga} - F_{\ga})^2 } \lesssim_\ga \frac{1}{(\al^2  n )^{(\ga - 1)\wedge 1}}
+ \mathbf{1}_{\{  \ga \geq 3\}}\frac{\| \ptau \|^{2(\ga - 2 ) }_{\ga -2}}{(\al^2  n )^2}
+ \mathbf{1}_{\{ \ga \geq 2\}}\frac{\| \ptau \|^{2\ga -2}_{2\ga -2}}{\al^2  n },
\end{equation*}
and thus is uniformly bounded over all $p\in \Pcal_K$ by
\begin{equation*}
     \sup_{p\in \Pcal_K} \E \cro{ (\widetilde{F}_{\ga} - F_{\ga})^2 } \lesssim_\ga
      \frac{1}{(\al^2n)^{\ga - 1}} + \frac{1}{\al^2n}.
\end{equation*}
\end{thm}

%The proof of the Theorem is written in Appendix. % \ref{appendic:proof:bias:interac}.

For $\ga >1$, the rate of the two-step procedure is therefore independent of the alphabet size $K$, unlike the rates of the plug-in and thresholded estimators. Hence, when $K$ is large, this rate is faster than those of the two (non-interactive) estimators. In particular, it is equal to $(\al^2 n)^{-1}$ as soon as $\ga \geq 2$, which is the minimax optimal rate (see next section).  Conversely for small $K$ and $\ga \in(1,2)$, this rate is slower than those of the plug-in and thresholded estimators. Accordingly, none of the three estimators is overall better than the others, and we discuss the choice of estimator in next section.

\subsection{Optimality of the Results}

Among the three estimators we proposed, the choice of estimators depends on the problem parameters $K$, $\ga$. The following recipe leads to a better estimator $\hat E_\gamma$. If $K\leq \sqrt{\al^2 n}$, define $\hat E_\gamma$ as the plug-in estimator $\hat F_\ga$;  otherwise (when $K> \sqrt{\al^2 n}$), $\hat E_\gamma$ is equal to the thresholded estimator $\overline{F}_\ga$ for $\ga <1$, and to the sequentially interactive procedure $\widetilde{F}_\ga$ for $\ga >1$. 
%In the next theorem,  we bound the maximal quadratic risk over the class $\mathcal{P}_K$ the estimator $\hat{E}_{\ga}$. 

\begin{cor}\label{mainTHM} The quadratic risk of $\hat{E}_{\ga}$  is uniformly bounded over all $p\in \mathcal{P}_K$ by \\
\textup{1)}  $\ga \in(0,1)$:
\begin{equation}\label{conj01}
    \sup_{p\in \Pcal_K} \E \cro{ (\hat{E}_{\ga} - F_{\ga})^2 } \lesssim_\ga K^{2(1-\ga )} \wedge \frac{K^2}{(\al^2 n)^{\ga }} \enspace.
\end{equation}
\textup{2)}  $\ga \in(1,3/2)$:
\begin{equation}\label{conj132}
    \sup_{p\in \Pcal_K} \E \cro{ (\hat{E}_{\ga} - F_{\ga})^2 } \lesssim_\ga  \frac{1}{(\al^2 n)^{\ga -1}} \wedge \left( \frac{K^2}{(\al^2 n)^{\ga }} + \frac{K^{3-2\ga}  }{\al^2  n } \right) \enspace.
\end{equation}
\textup{3)}  $\ga \in(3/2,2)$:
\begin{equation}\label{conj322}
     \sup_{p\in \Pcal_K}  \E \cro{ (\hat{E}_{\ga} - F_{\ga})^2 } \lesssim_\ga  \frac{1}{(\al^2 n)^{\ga -1}} \wedge \left( \frac{K^2}{(\al^2 n)^{\ga }} + \frac{1  }{\al^2  n } \right) \enspace.
\end{equation}
\textup{4)}  $\ga \geq 2$:
\begin{equation}\label{conj2}
    \sup_{p\in \Pcal_K} \E \cro{ (\hat{E}_{\ga} - F_{\ga})^2 } \lesssim_\ga  \frac{1}{\al^2  n}  \enspace.
\end{equation}
\end{cor}

The rate  \eqref{conj2} for $\ga \geq 2$ can be seen as the private parametric rate, which is minimax optimal (see lower bounds in Theorem \ref{simple:LB:thm}). Then for $\ga \in(0,2)$, the smaller $\ga$ is, the slower the rates are. Specifically, each upper bound in (\ref{conj01}) to  (\ref{conj322}) is the minimum of two error bounds, with a phase transition at $K \asymp \sqrt{\al^2 n}$. Above this transition level, i.e. when $K\gtrsim \sqrt{\al^2 n}$, the upper bounds (\ref{conj132}) and (\ref{conj322}) for $\ga \in (1,2)$ are equal to the first error bound $(\al^2 n)^{-(\ga -1)}$ which is free of $K$. Below the transition level, i.e. when  $K\lesssim \sqrt{\al^2 n}$, they are equal to the second error bound that depends on $K$. Hence, for $\ga \in (1,2)$, the upper bounds on the maximal quadratic risk of $\hat{E}_{\ga}$ increase with $K$ as long as $ K \lesssim \sqrt{\al^2 n}$, and become equal to $(\al^2 n)^{-(\ga -1)}$ for any $K$ larger than  $\sqrt{\al^2 n}$. In contrast, the upper bound \eqref{conj01} for $\ga \in(0,1)$ depends on $K$ regardless of the value of  $K$. Thus, the rates get slower as the power $\ga$ decreases, where $\ga$ can be seen as a smoothness indicator of the function $x \mapsto x^{\ga}$ and thus of the functional $F_{\ga} =\sum_k p_k^{\ga}$.

%We outline the rest of the paper where three estimators will be presented. The two first estimators use a non-interactive PM which is introduced in section 3.1. The first estimator we study  is a plug-in estimator defined in section 3.2. Its maximal quadratic risk is bounded from above by  the second error bound in (\ref{conj01}), (\ref{conj132}), (\ref{conj322}), and thus suffers from a  quadratic dependence on $K$ for all $K$, which will turn out to be a tight characterization of the performance of the plug-in estimator (see Proposition~\ref{borne:inf:plugin:petit:pk}). In order to reduce this dependence on $K$, we will build on this first estimator and define a thresholded estimator in section 3.3. The quadratic risk of this  thresholded estimator matches \eqref{conj01} exactly, and (\ref{conj132}), (\ref{conj322}), (\ref{conj2}) only up to $\log$ factors. We will  remove this $\log$ factor in section 4 by using a sequentially interactive PM, instead of the  non-interactive PM used before (in section 3). However, this sequentially interactive procedure  only achieves the rates (\ref{conj132}-\ref{conj322}) for large $K$ (s.t. $K\geq \sqrt{\al^2 n}$). See Table~\ref{Rates} for a summary. Hence, none of the three estimators above achieves alone all the desired rates (\ref{conj01}) to (\ref{conj2}), and accordingly . We give  lower bounds with respect to all estimators and all PM in section 5. 

%%%%%%%%%%%%%%%%%%%%%%%
%%%%%%%%%%%%%%%%%%%%%%

We give lower bounds over all estimators and all $\alpha$-LDP   sequentially interactive  mechanisms. Recall that non-interactive mechanisms are just a special case of sequentially interactive mechanisms, thus our lower bounds hold in particular for non-interactive mechanisms. In the special case where $K \lesssim 1$ is a numerical constant, the next lower bounds match the rate of the plug-in estimator (Theorem \ref{thm:bias:nonInterac:summary:new}) for any $\ga >0$, $\ga \neq 1$.

\begin{thm}\label{simple:LB:thm}For any $\ga >0$, $\ga\neq 1$ and integer $K\geq 2$, we have the lower bound
\begin{equation*}
    \underset{Q, \hat{F}}{\textup{inf}} \ \, \underset{p  \in \Pcal_K}{\textup{sup}} \, \E \cro{ ( \hat{F} - F_{\ga})^2 }  \gtrsim_\ga \frac{1}{\al^2 n} + \frac{1}{(\al^2 n)^{\ga}} \enspace,
\end{equation*}
where the infimum is taken over all estimators $\hat{F}$ and all $\al$-LDP sequentially interactive PM $Q$.
\end{thm}

In the general situation of any $K \geq 2$, we thus have the minimax optimal rate $(\al^2 n)^{-1}$ when $\ga \geq 2$, which is achieved by the sequentially interactive estimator (Theorem \ref{thm:bias:summary:interac}).   In contrast for $\gamma \in (0,2)$, the next lower bounds  depend on the alphabet size $K$. 

\begin{thm}\label{nosimple:LB:thm} %$\al > 0$.
For any integer $K \geq 2$, we have the lower bounds\\
\textup{1)} $\ga\in(0,1)$:
\begin{equation}\label{THM:eq:LB}
    \underset{Q, \hat{F}}{\textup{inf}} \ \,  \underset{p  \in \Pcal_K}{\textup{sup}} \, \E \cro{ ( \hat{F} - F_{\ga})^2 }  \gtrsim_\ga  K^{2(1-\ga)}  \wedge  \frac{K^{2-\ga}}{ ( n (e^{2\al}-e^{-2\al})^2    )^{\ga}}
\end{equation}
\textup{2)} $\ga\in(1,2)$:
\begin{equation}\label{THM:eq:LB:geq1}
    \underset{Q}{\textup{inf}} \, \underset{\hat{F}}{\textup{inf}} \, \underset{p  \in \Pcal_K}{\textup{sup}} \, \E \cro{ ( \hat{F} - F_{\ga})^2 }  \gtrsim_\ga  \frac{1}{ [  (e^{2\al}-e^{-2\al})^2  n  ]^{2(\ga-1)}} \wedge  \frac{K^{2-\ga}}{ [  (e^{2\al}-e^{-2\al})^2   n ]^{\ga}}
\end{equation}
where the infimum is taken over all estimators $\hat{F}$ and all $\al$-LDP sequentially interactive PM $Q$.
\end{thm}

Note that $e^{2\al}-e^{-2\al} \asymp \al$ when $\al$ is small.
Then the lower bound \eqref{THM:eq:LB} for $\ga \in(0,1)$ matches, up to a factor $K^{\ga}$, the quadratic risk of the thresholded estimator (Theorem \ref{threshold:thm}), or equivalently the rate \eqref{conj01} of the combined estimator $\hat{E}_{\ga}$.  For $\ga \in(1,2)$, the situation is slightly more involved, but one can  see again that there is a gap of a factor $K^{\ga}$  between some terms of the lower bound \eqref{THM:eq:LB:geq1} and the upper bounds (\ref{conj132}) and (\ref{conj322}) of Corollary \ref{mainTHM}. Besides, reading the proof of lower bounds, one can check that this $K^{\ga}$-gap  is tightly connected to the gap  between the left-terms $1/(\al^2n)^{2(\ga-1)}$ and $1/(\al^2n)^{\ga-1}$ of the lower and upper bounds respectively, so that the optimality of the problem actually boils down to reducing this  $K^{\ga}$-gap.

\subsection{Discussion} In this first work on the estimation of the power sum functional $F_{\ga}$, we focus on a plug-in estimator based on samples obtained using a non-interactive PM. Although its analogue estimator in the non private case is known to be nearly minimax optimal, we show that our plug-in estimator performs poorly in the LDP setting for large $K$ (proving a tight characterization of its maximal quadratic risk). Hence, practical methods performing well in the non private case should not be systematically transferred to the LDP setting, which requires to design new statistical procedures.  We then suggest a correction of this estimator using thresholding, which significantly improves the rates of convergence, removing almost  the whole dependence of the risk in the support size $K$ for large $K$. We finally get faster rates by combining these two previous estimators with a sequentially interactive procedure. It is also important to highlight that all privacy mechanisms and estimators introduced here could be further investigated for other functionals, for example the Rényi entropy $H_{\ga}$ which is  fundamental in information theory and is connected to the power sum $F_{\ga}$ via the relation $H_{\ga}= \frac{\log F_{\ga}}{1-\ga}$.

We conjecture that our upper bounds are tight, up to some logarithmic factors. In future work, this logarithmic factor should be reduced via the best polynomial approximation method. This is similar to the line of work in the standard setup (without privacy constraint), where the analogue of our plug-in estimator (the MLE) is known to be nearly minimax optimal up to a poly-logarithmic factor [14], while the the best polynomial approximation method closes this logarithmic-gap and is minimax optimal \cite{Jiao_2015,WuYang}. 
Leaving aside the logarithmic factors, we conjecture that our lower bounds are not optimal, up to a factor $K^{\ga}$. The two fuzzy hypothesis theorem is often used in the standard (non-private) setting to derive lower bounds on the estimation rate of functionals such as the power sum  $F_{\ga}$. A challenge in the LDP setting is to provide such turnkey tools that help prove universal lower bounds.   

Another challenge is to understand when sequentially interactive procedures outperform non-interactive ones.  When $\ga >1$, we have gained logarithmic factors in our error bounds  by considering sequentially interactive procedures. However, because the optimal rate of non-interactive mechanisms is not proven, it is unclear that our logarithmic gap between non-interactive and sequentially interactive actually exists.
The estimation of power sum functionals in the context of local differential privacy proves to be a rich topic potentially difficult to solve sharply in all possible cases.

\bibliographystyle{plain}

%%%%%%%%%%%%%%%%%%%%%%%%%%%%%%%%%%%%%%%%%%%%%%%%%%%%%%%%%%%%

\newpage

%%%%%%%%%%%%%%%%%%%%%%%%%%%%%%%%%%%%%%%%%%%%%%%%%%%%%%%%%%%%

\centerline{{\huge{\textbf{Supplementary Material}}}}

\bigskip 

\bigskip 

\bigskip

This Supplementary Material contains the proofs of the results presented above.

\appendix

\section{Proofs of upper bounds }

%%%%%%%%%%%%%%%%%%%%%%%%%%%%%%%%
\subsection{Plug-in estimator }\label{proof:sect:thm}

\noindent \textbf{Proof of $1^{\text{st}}$ bound in  Theorem 2.1.} \textit{$1^{\circ}$. Bias:} We have using the triangle inequality,
\begin{equation*}
    \left|  \E \hat{F}_{\ga} - F_{\ga} \right|=   \left|  \E \hat{F}_{\ga} - \sum_{k=1}^K p_k^{\ga} \right| \leq \sum_{k=1}^K \left| \E \hat{F}_{\ga}(k) - p_k^{\ga} \right| \enspace.
\end{equation*}
Hence, it suffices to upper bound the $k^{\textup{th}}$ bias component $|\E \hat{F}_{\ga}(k) - p_k^{\ga}|$ for all $k\in [K]$ and $\ga \neq 1$ (the case $\ga =1$ being trivial). We separate the analysis in two different ranges of values of $p_k$. Define $\Kcal_{\geq \tau} = \{ k\in [K]:\, p_k \geq \tau \}$, and $\Kcal_{< \tau} = [K] \setminus \Kcal_{\geq \tau} $. By Lemma \ref{lem:bias:small:pk} we have
\begin{equation*}
    \sum_{k\in \Kcal_{< \tau} } \left| \E \hat{F}_{\ga}(k) - p_k^{\ga} \right| \leq  C \frac{ |\Kcal_{< \tau}| }{(\al^2 n)^{\ga/2}}
\end{equation*}for a constant $C$ depending only on $\ga$.
Lemma \ref{lem:bias:big:pk} ensures that
\begin{equation*}
    \sum_{k\in \Kcal_{\geq \tau} } \left| \E \hat{F}_{\ga}(k) - p_k^{\ga} \right| \leq  C' \pa{  \frac{| \Kcal_{\geq \tau}|}{(\al^2  n)^{\ga/2}}  + \mathbbm{1}_{\{  \ga \geq 2\}}\frac{\| \ptau \|^{\ga -2}_{\ga -2}}{\al^2  n }  }
\end{equation*}for a constant $C'$ depending only on $\ga$. Gathering the above inequalities, we have
\begin{equation}\label{proof:bias:eq:fina}
      \left|  \E \hat{F}_{\ga} - F_{\ga} \right| \leq  (C+C') \pa{  \frac{K}{(\al^2 n)^{\ga/2}}  + \mathbbm{1}_{\{  \ga \geq 2\}}\frac{\| \ptau \|^{\ga -2}_{\ga -2}}{\al^2  n }  } .
\end{equation}

%%%%%%%%%

\textit{$2^{\circ}$. Variance:} By Lemma \ref{lem:negassociattion} , we have $ \textup{Cov}\big{(}\hat{F}_{\ga}(k), \hat{F}_{\ga}(k') \big{)} \leq 0$ for any $k \neq k' \in[K]$. Hence
\begin{equation}\label{VarDecomposition:NegativeAssociation}
    \Va \pa{ \sum_{k=1}^K  \hat{F}_{\ga}(k) } \leq \sum_{k=1}^K  \Va \pa{  \hat{F}_{\ga}(k) } \enspace.
\end{equation}

As in the  proof of the bias bound above, we separate our analysis in two different ranges of values of $p_k$. For small $p_k$, we use Lemma \ref{lem:bias:small:pk} to get
\begin{equation*}
    \sum_{k\in \Kcal_{< \tau} }  \Va \pa{  \hat{F}_{\ga}(k) }   \leq  \widetilde C \frac{| \Kcal_{< \tau}|}{(\al^2 n)^{\ga}} \enspace,
\end{equation*}where $\widetilde C$ is a constant depending only on $\ga$. For large $p_k$, we deduce from Lemma \ref{lem:var:big:pk} that
\begin{equation*}
      \sum_{k\in \Kcal_{\geq \tau} } \Va \pa{  \hat{F}_{\ga}(k) } \leq   \widetilde C' \bigg{(} \frac{|\Kcal_{\geq \tau}|}{(\al^2 n)^{\ga}}  + \mathbbm{1}_{\{  \ga \geq 1\}}\frac{\| \ptau \|^{2\ga -2}_{2\ga -2}}{\al^2  n } \bigg{)}
\end{equation*}
for a constant $\widetilde C'$ depending only on $\ga$. Then, plugging these bounds into \eqref{VarDecomposition:NegativeAssociation}, we have
\begin{equation}\label{proof:var:eq:fina}
    \Va \pa{ \sum_{k=1}^K  \hat{F}_{\ga}(k) } \leq  (\widetilde C + \widetilde C') \bigg{(} \frac{K}{(\al^2  n)^{\ga}}  + \mathbbm{1}_{\{  \ga \geq 1\}}\frac{\| \ptau \|^{2\ga -2}_{2\ga -2}}{\al^2  n } \bigg{)} \enspace.
\end{equation}
%%%%%%%%%%
The proof of the of $1^{\text{st}}$ bound in Theorem 2.1 is complete. \hfill $\square$

 %%%%%%%%%%%%%
 \bigskip
 %%%%%%%%%%%%

 \noindent \textbf{Proof of $2^{\text{nd}}$ bound in  Theorem 2.1.} We only need to control the second and third terms of the  $1^{\text{st}}$ bound in Theorem 2.1. The squared root of the second term is bounded from above by 
\begin{equation*}
     \frac{ \sum_{k=1}^K p_k^{\ga -2} \mathbbm{1}_{\{p_k > \tau\}}}{\al^2  n }
     \leq \sum_{k=1}^K \frac{  p_k^{\ga -1} }{\sqrt{\al^2  n} }  \frac{ p_k^{-1} \mathbbm{1}_{\{p_k > \tau\}} }{\sqrt{\al^2 n} }
     \leq \sum_{k=1}^K \frac{  p_k^{\ga -1}  c^{-1}}{\sqrt{\al^2  n} }
     =  \frac{  \| p \|^{\ga -1}_{\ga -1} c^{-1} }{\sqrt{\al^2  n} } \enspace.
\end{equation*}
Since $(p_k)_k$ are probabilities, we have $p_k^{\ga - 1} \leq p_k$ for $\ga \geq 2$ and we can further bound the last display by $ \| p \|^{\ga -1}_{\ga -1} \leq \sum_{k=1}^K p_k = 1$ for $\ga \geq 2$. Hence, the second term is bounded by $\mathbbm{1}_{\ga \geq 2}(\al^2n)^{-1}$.

Let us bound the third term. Since $\sum_k p_k =1$, the number of the significant $p_k \geq \tau$ is necessarily smaller than $\tau^{-1} = c^{-1} \sqrt{\al^2 n}$, and thus smaller than $K_{\wedge \tau^{-1}} := K \wedge \sqrt{\al^2  n}$. Then, when $\ga \in (1, 3/2)$,  we use the concavity to have  $\| \ptau \|^{2\ga -2}_{2\ga -2}  \leq K_{\wedge \tau^{-1}}^{3 -2\ga}$ for all $p\in \Pcal_K$. When $\ga \geq 3/2$ we have  $\| \ptau \|^{2\ga -2}_{2\ga -2} \leq 1$. Therefore, the third term is uniformly bounded over the class $\Pcal_K$ by
\begin{equation*}
   \mathbbm{1}_{\{  \ga \geq 1\}}\frac{\| \ptau \|^{2\ga -2}_{2\ga -2}}{\al^2  n }  \leq  \mathbbm{1}_{\{  \ga \geq 1\}}\frac{1 \vee K_{\wedge \tau^{-1}}^{3-2\ga}}{\al^2  n }  \enspace.
\end{equation*}

This concludes the proof of the $2^{\text{nd}}$ bound in  Theorem 2.1. \hfill $\square$

%%%%%%%%%%%%%%%%

 \subsection{Thresholded plug-in estimator (proof of Theorem 2.3)}\label{section:proof:threshold}

\underline{Case $\ga \in (0,1)$}: Let us check the first bound of Theorem 2.3. We use the concavity of the power function $p^\gamma$ to have  $F_\ga \leq K (\sum_{k=1}^K p_k / K)^\ga = K^{1-\ga}$. Then, the quadratic risk of the trivial estimator $0$ is bounded by $K^{2(1-\ga)}$. On the other hand, the quadratic risk of the plug-in $\hat{F}_{\ga}$ is bounded by $K^2/(\al^2 n)^{\ga}$ (Theorem 2.1). Therefore,  the quadratic risk of the thresholded estimator $\overline{F}_{\ga} := \mathbbm{1}_{K\leq \tau^{-1}} \hat{F}_{\ga}$ satisfies the first bound of Theorem 2.3.

\underline{Case $\ga > 1$}: Recall that  $\hat{\tau} \asymp \sqrt{\log(Kn)/(\al^2  n)}$. We will prove the next bound on the risk of $\overline{F}_{\ga}$,
\begin{equation}\label{proof:eq:threshold:interm}
  \E \cro{ (\overline{F}_{\ga} - F_{\ga})^2 } 
      \lesssim_{\ga}  (K\hat{\tau}^{\ga} \wedge \hat{\tau}^{\ga -1})^2 + \frac{ (K \wedge \hat{\tau}^{-1})^{3-2\ga} \vee 1 }{\al^2  n } \enspace.
\end{equation}Before that, we  check that \eqref{proof:eq:threshold:interm} implies the second inequality of Theorem 2.3. 

\textit{(i) Assume that  $K \geq \hat{\tau}^{-1}$}, then the RHS of \eqref{proof:eq:threshold:interm} becomes 
\begin{equation*}
     \hat{\tau}^{2(\ga -1)} + \frac{ \hat{\tau}^{2\ga-3} \vee 1 }{\al^2  n } \lesssim    \frac{\left(\log(Kn)\right)^{\ga-1}}{(\al^2 n)^{\ga-1}} +\frac{\left(\log(Kn)\right)^{\ga-(3/2)} }{ (\al^2 n)^{\ga-(1/2)}  } + \frac{1 }{\al^2  n }  \lesssim    \frac{\left(\log(Kn)\right)^{\ga-1}}{(\al^2 n)^{\ga-1}} + \frac{1 }{\al^2  n }  \enspace,
\end{equation*}where the last inequality follows from the bound
\begin{equation*}
       \frac{\left(\log(Kn)\right)^{\ga-(3/2)} }{ (\al^2 n)^{\ga-(1/2)}  } \leq    \frac{\left(\log(Kn)\right)^{\ga-1}}{(\al^2 n)^{\ga-1}} \enspace,
\end{equation*}which is equivalent to $ \al^2 n \log(Kn)  \geq  1 $. Hence, \eqref{proof:eq:threshold:interm} is upper bounded by the smallest term of the second inequality of Theorem 2.3.

\textit{(ii) Assume that  $K \leq \hat{\tau}^{-1}$}, then  the RHS of \eqref{proof:eq:threshold:interm} becomes 
\begin{equation*}
     K^2 \hat{\tau}^{2\ga} + \frac{ K^{3-2\ga} \vee 1 }{\al^2  n } \lesssim    \frac{K^2 \left(\log(Kn)\right)^{\ga}}{(\al^2 n)^{\ga}} +  \frac{1 \vee K^{3-2\ga}  }{\al^2  n } \enspace,
\end{equation*}which is the smallest term of the second inequality of Theorem 2.3. Hence, we have proved that the second inequality of Theorem 2.3 follows from  \eqref{proof:eq:threshold:interm}.

 \medskip
 %%%%%%%%%%%%

 \noindent \textbf{Proof of \eqref{proof:eq:threshold:interm}.}  We have the deterministic bound
\begin{equation*}
    |\ovf_{\ga} - F_{\ga}| \leq \ovf_{\ga} + F_{\ga} \leq K(2^{\ga}+1) \enspace.
\end{equation*}
Introduce the following event
\begin{equation*}
    A= \left\{ \exists k\in[K]:\, \left(\hat{z}_k^{(1)} < \hat{\tau} \textup{ and } p_k  \geq 3\hat{\tau}/2 \right) \, \textup{or} \, \left(\hat{z}_k^{(1)} \geq \hat{\tau} \textup{ and } p_k  < \hat{\tau}/2  \right) \right\}
\end{equation*} and denote  the complementary event by $A^c$. We have
\begin{equation}\label{proof:threslhold:eq1}
    \E \cro{(\ovf_{\ga} - F_{\ga})^2} \leq  \E \cro{\mathbbm{1}_{A^c}(\ovf_{\ga} - F_{\ga})^2} + \Pb(A)\left(K(2^{\ga}+1)\right)^2 \enspace.
\end{equation}
Let us bound the second term of the RHS of \eqref{proof:threslhold:eq1} by showing that $\Pb(A)\leq 6/(K^2n)$.  By assumption in the theorem, we have  $n \geq 2 \log(K)$. This ensures that $n \geq \log(K n^{1/3})$, which allows us to use Lemma \ref{hoeffding:zk} which gives  $\Pb \left(|\hat{z}_k^{(1)} - p_k | > \hat{\tau}/2 \right)  \leq6/(K^3n) $. Hence, for $p_k \geq 3\hat{\tau}/2 $, we have
\begin{equation*}
      \Pb \left(\hat{z}_k^{(1)} <  \hat{\tau}\right)  \leq \frac{6}{K^3n}   \enspace,
\end{equation*}
and for $p_k < \hat{\tau}/2$,
\begin{equation*}
      \Pb \left(\hat{z}_k^{(1)} \geq  \hat{\tau}\right)  \leq \frac{6}{K^3n}   \enspace .
\end{equation*}
We then use the union bound over $k\in[K]$ to get $\Pb(A) \leq 6/(K^2n)$. The second term of the RHS of \eqref{proof:threslhold:eq1} is therefore bounded by $6(2^{\ga}+1)^2/n$.

We now control the first term of the RHS of \eqref{proof:threslhold:eq1}.
For any real $a>0$ ,  we note $\Kcal_{< a} = \{ k\in [K]: \, p_k < a \}$ and $\hat{\Kcal}_{< a} = \{ k\in [K]: \, \hat{z}_k^{(1)} < a \}$, with  their respective complementary sets $\Kcal_{\geq a} = [K] \setminus \Kcal_{< a}$ and $\hat{\Kcal}_{\geq a} =  [K] \setminus \hat{\Kcal}_{< a}$. Splitting the sum over the $k$ in $\hat{\Kcal}_{ < \hat{\tau}} $ and  $\hat{\Kcal}_{ \geq \hat{\tau}} $ respectively, we get
\begin{equation*}
    \mathbbm{1}_{A^c}(\ovf_{\ga} - F_{\ga})^2 \leq 2  \mathbbm{1}_{A^c} \pa{\|(p_k)_{k \in \hat{\Kcal}_{ < \hat{\tau}} }\|_{\ga}^{\ga}}^2 +  2 \mathbbm{1}_{A^c}\Big{(}\sum_{k\in \hat{\Kcal}_{\geq \hat{\tau}} }\ovf_{\ga}(k) - F_{\ga}(k)\Big{)}^2 \enspace.
\end{equation*}Since  $\hat{\Kcal}_{< \hat{\tau}} \subset \Kcal_{< 3\hat{\tau}/2}$ on the event $A^c$, we can bound the first term  by 
\begin{equation*}
   \mathbbm{1}_{A^c} \|(p_k)_{k \in \hat{\Kcal}_{ < \hat{\tau}} }\|_{\ga}^{\ga} \leq  \|(p_k)_{k \in \Kcal_{ <3\hat{\tau}/2} }\|_{\ga}^{\ga} \leq K(3\hat{\tau}/2)^{\ga} \wedge (3\hat{\tau}/2)^{\ga-1}
\end{equation*} for any $\ga >1$ and $p\in \mathcal{P}_K$. For the second  term,  we will use the independence between the data samples $z^{(1)}:= (z^{(1)}_1,\ldots,z_n^{(1)})$  and $z^{(2)}:= (z^{(2)}_1,\ldots,z_n^{(2)})$. In particular, the set $\hat{\Kcal}_{\geq \hat{\tau}}$ and the event $A^c$ are deterministic conditionally to $z^{(1)}$, so that
\begin{align*}
  \E \cro{\mathbbm{1}_{A^c}\Big{(}\sum_{k\in \hat{\Kcal}_{\geq \hat{\tau}} }\ovf_{\ga}(k) - F_{\ga}(k)\Big{)}^2 \Big{|} z^{(1)} } 
  &= \mathbbm{1}_{A^c} \E \cro{\Big{(}\sum_{k\in \hat{\Kcal}_{\geq \hat{\tau}} }\ovf_{\ga}(k) - F_{\ga}(k)\Big{)}^2 \Big{|} z^{(1)} }\\
  &\leq  \mathbbm{1}_{A^c} C \pa{  \frac{|\hat{\Kcal}_{\geq \hat{\tau}}|^2}{(\al^2n)^{\ga}}   + \frac{ |\hat{\Kcal}_{\geq \hat{\tau}}|^{3-2\ga} \vee 1}{\al^2  n } }
\end{align*}where the last line is similar to the $2^{\text{nd}}$ bound in Theorem 2.1  with $K$ replaced by $|\hat{\Kcal}_{\geq \hat{\tau}}|$, and where $C$ is some constant depending only on $\ga$. We can further bound the last display by noting that $\hat{\Kcal}_{\geq \hat{\tau}} \subset \Kcal_{\geq \hat{\tau}/2}$ on the event $A^c$, and $|\Kcal_{\geq \hat{\tau}/2}| \leq K \wedge ( \hat{\tau}/2)^{-1}$. Going back to \eqref{proof:threslhold:eq1}, we then have for all $p\in \Pcal_K$,
\begin{align*}
  \E \cro{ (\overline{F}_{\ga} - F_{\ga})^2 } &\lesssim_{\ga}  (K\hat{\tau}^{\ga} \wedge \hat{\tau}^{\ga -1})^2
      + \frac{(K \wedge \hat{\tau}^{ -1})^2}{(\al^2  n)^{\ga }}
      +  \frac{ (K \wedge \hat{\tau}^{-1})^{3-2\ga} \vee 1 }{\al^2  n } + \frac{1 }{ n }\\
      &\lesssim_{\ga}  (K\hat{\tau}^{\ga} \wedge \hat{\tau}^{\ga -1})^2 + \frac{ (K \wedge \hat{\tau}^{-1})^{3-2\ga} \vee 1  }{\al^2  n } \enspace.
\end{align*} The proof of \eqref{proof:eq:threshold:interm} is complete. \hfill $\square$

\subsection{Interactive privacy mechanism }\label{appendic:proof:bias:interac}

\textbf{Proof of $1^{\text{st}}$ bound in Theorem 2.4.} \textit{$1^{\circ}$. Bias:} We decompose the expected value of $\widetilde{F}_{\ga}$ :
\begin{align}\label{proof:bias:decomp}
    \E \widetilde{F}_{\ga} &=  \frac{1}{n} \sum_{i=1}^n\E \E \cro{ z_i^{(2)} | z^{(1)}, z^{(2)}}  =   \frac{1}{n} \sum_{i=1}^n\E \E \cro{\hat{F}_{\ga-1}^{(1)}(x_i^{(2)}) | z^{(1)}, x^{(2)}} \nonumber \\
   &= \sum_{k=1}^K p_k \E \E \cro{\hat{F}^{(1)}_{\ga-1}(k)| z^{(1)}} = \sum_{k=1}^K p_k \E \cro{\hat{F}^{(1)}_{\ga-1}(k)}
\end{align}
so that, for any $\ga >1$, $\ga \neq 2$ (the case $\ga =2$ being trivial), we have
\begin{align}\label{proof:interBias}
    \left|  \E \widetilde{F}_{\ga} - \sum_{k=1}^K p_k^{\ga} \right| &\leq \sum_{k=1}^K p_k \left| \E \hat{F}^{(1)}_{\ga-1}(k) - p_k^{\ga-1} \right| \nonumber \\
    &\leq C \pa{  \frac{1}{(\al^2  n)^{(\ga-1)/2}}  + \mathbbm{1}_{\{  \ga \geq 3\}}\frac{\| \ptau \|^{\ga -2}_{\ga -2}}{\al^2 n }  }
\end{align}using Lemma \ref{lem:bias:small:pk} and \ref{lem:bias:big:pk} and $\sum_k p_k =1$, where $C$ is a constant depending only on $\ga$.

\textit{$2^{\circ}$. Variance:} By the law of total variance we have
\begin{equation}\label{Variance:decomposition:LawTotalVariance}
    \textup{Var}\pa{\widetilde{F}_{\ga}} = \E \cro{\textup{Var} \pa{\widetilde{F}_{\ga} |z^{(1)}}} + \Va \pa{\E\cro{\widetilde{F}_{\ga}|z^{(1)}}} \enspace.
\end{equation}
We control the first term in the RHS of \eqref{Variance:decomposition:LawTotalVariance}:
\begin{align*} 
    \Var  \pa{\widetilde{F}_{\ga} |z^{(1)}} &= \frac{1}{n} \Var\pa{ z_1^{(2)} |z^{(1)}} \leq \frac{1}{n} \E\cro{ \left(z_1^{(2)}\right)^2 |z^{(1)}} \nonumber \\
    &= \frac{2^{2\ga -1}}{n}  \pa{\frac{e^{\al}+1}{e^{\al}-1}}^2 \leq  \frac{2^{2\ga +1}}{\al^2 n}
\end{align*}
where we used
$(\frac{e^{\al}+1}{e^{\al}-1})^2 = (1 + \frac{1}{e^{\al}-1})^2 \leq (1 + \frac{1}{\al})^2 \leq \frac{4}{\al^2}$.
For the second term in the RHS of \eqref{Variance:decomposition:LawTotalVariance}, we have using \eqref{proof:bias:decomp}
\begin{equation*}
    \Va \pa{\E\cro{\widetilde{F}_{\ga}|Z^{(1)}}} = \Va \pa{ \sum_{k=1}^K p_k  \hat{F}_{\ga-1}^{(1)}(k) }  \leq \sum_{k=1}^K p_k^2  \Va \pa{ \ \hat{F}_{\ga-1}^{(1)}(k) } 
\end{equation*}
where the inequality can be deduced from Lemma \ref{lem:negassociattion}. Then, by Lemma \ref{lem:bias:small:pk} and \ref{lem:var:big:pk},
\begin{equation*}
   \sum_{k=1}^K p_k^2  \Va \pa{ \ \hat{F}_{\ga-1}^{(1)}(k) } \leq   \widetilde C \bigg{(} \frac{\| p \|^{2}_{2}}{(\al^2  n)^{\ga -1}}  + \mathbbm{1}_{\{  \ga \geq 2\}}\frac{\| \ptau \|^{2\ga -2}_{2\ga -2}}{\al^2 n } \bigg{)}  
\end{equation*}
for a constant $\widetilde C $ depending only $\ga$. The proof of the $1^{\text{st}}$ bound in  Theorem 2.4 is complete. \hfill $\square$

%%%%%%%%%%%%%%

\bigskip

\noindent \textbf{Proof of $2^{\text{nd}}$ bound in Theorem 2.4.}  The desired bound  follows from the $1^{\text{st}}$ bound of  Theorem 2.4 and the fact that  $\mathbbm{1}_{\{  \ga \geq 3\}} \| \ptau \|^{\ga -2}_{\ga -2} \leq 1$ and $\mathbbm{1}_{\{  \ga \geq 2\}} \| \ptau \|^{2\ga -2}_{2\ga -2} \leq \| p \|^{2}_{2} \leq 1$ for all $p\in \Pcal_K$.   \hfill $\square$

 %%%%%%%%%%%

%%%%%%%%%%%%%%%%%%%

\section{Main lemmas for upper bounds} \label{appen:section:lemm}
We use the notations $\ovx_k=  \frac{1}{n} \sum_{i=1}^n \mathbbm{1}_{\{x_i = k\}}$ and $\ovw_k =  \frac{1}{n} \sum_{i=1}^n w_{ik}$, so that $\hat{z}_k =\ovx_k + \frac{\si}{\al} \ovw_k$. We consider $\al \in(0,\infty)$  in this Appendix \ref{appen:section:lemm},  unlike in the main section of the paper where we assumed that $\al\in(0,1)$ and $\al^2n \geq 1$.

\subsection{Concentration of $\hat{z}_k$}

We control  the  concentration of $\hat{z}_k$ in the next lemma.

\begin{lem}\label{concentrationBinomialEtLaplace}
For any $\al \in (0,\infty)$ and any $r>0$, we have
\begin{align*}
    \E\cro{ |\hat{z}_k-p_k|^r} &\leq \frac{C_{BL,r}}{((\al^2 \wedge 1)n)^{r/2}} \enspace,\\
     \E \cro{ |\hat z_{k}|^{r}} &\leq \frac{2^{r} C_{BL,r}}{((\al^2 \wedge 1)n)^{r/2}} +  2^{r} p_k^{r} \enspace,
\end{align*}where $C_{BL,r}$ is a constant depending only on $r$. Besides,
\begin{align*}
   \Pb (\hat{z}_k < \frac{p_k}{2})  &\leq 3 \, \exp\cro{-\frac{n}{128 } \bigg{(}\frac{(\al \wedge 1) p_k}{\si}\bigg{)}^2  }   \enspace.
\end{align*}
\end{lem}

\medskip

\noindent \textbf{Proof of Lemma \ref{concentrationBinomialEtLaplace}.}  By \eqref{ep:Mean:Binom} in Lemma \ref{hatZ:bigP:nonInterac} and \eqref{eq2:lemmaC2} in Lemma \ref{LaplaceControl}, we have for any $r>0$,
\begin{align*}
\E\cro{ |\hat{z}_k-p_k| ^r} \leq  2^{r} \E\cro{ |\ovx_k-p_k|^r} + 2^{r}  \E \cro{ \Big{(}\frac{ \si |\ovw_k|}{\al}\Big{)}^{r}}
   &\leq \frac{2^{r}C_{B,r}}{n^{r/2}}  +   \frac{(2\si)^r C_{L,r}}{(\al^2 n)^{r/2}} \\
   &\leq \frac{2^{r} \left(C_{B,r} + \si^r C_{L,r}\right)}{((\al^2 \wedge 1)n)^{r/2}}
\end{align*}where $C_{B,r}$ and $C_{L,r}$ are constants that only depend on $r$. Then, denoting $C_{BL,r}= 2^{r} \left(C_{B,r} + \si^r C_{L,r}\right)$, we have
\begin{align*}
   \E \cro{ |\hat z_{k}|^{r}} =  \E \cro{ |\hat z_{k}-p_k + p_k|^{r}} &\leq  2^{r}\E \cro{ |\hat z_{k}-p_k|^{r}} + 2^{r} p_k^{r}\\
  &\leq \frac{2^r C_{BL,r}}{((\al^2 \wedge 1)n)^{r/2}} + 2^{r} p_k^{r} \enspace.
\end{align*}
Finally, by  \eqref{concentratoinBinomialJiao} in Lemma \ref{hatZ:bigP:nonInterac} and \eqref{eq1:lemmaC2} in Lemma \ref{LaplaceControl}, we have
\begin{align*}
\Pb (\hat{z}_k < \frac{p_k}{2})  \leq\Pb(\ovx_k < \frac{3 p_k}{4  }) + \Pb(\frac{\si \ovw_k}{\al }< -\frac{ p_k}{4  })
    &\leq  e^{-\pa{\frac{1}{4}}^2\frac{np_k}{2}} + e^{-\frac{n}{8} \pa{\frac{\al p_k}{4\si }}^2 }   + e^{-\frac{n}{4} \pa{\frac{\al p_k}{4\si }} }  \\
     &\leq  3\,  e^{-\frac{n}{128 \si^2} \pa{(\al \wedge 1) p_k}^2 }  \enspace.
\end{align*}
The proof of Lemma \ref{concentrationBinomialEtLaplace} is complete. \hfill $\square$

\bigskip

Recall that $\hat{F}_{\ga}(k)= \left(T_{[0,2]}[\hat{z}_k]\right)^{\ga}$. We bound the difference between the expectations of $T_{[0,2]}[\hat{z}_k]$ and $\hat{z}_k$ in the next lemma.

\medskip

\begin{lem}\label{lem:control:first:order}
We have for any $\al \in (0,\infty)$,
\begin{align*}
     \left|\E \cro{T_{[0,2]}[\hat{z}_k]}-p_k \right| \leq \frac{ 2p_k^{-1}}{(\al^2 \wedge 1)n} \pa{ \si^2 C_{L,2} +  \frac{16\ga}{e}} \enspace.
\end{align*}
\end{lem}

\medskip

\noindent \textbf{Proof of Lemma \ref{lem:control:first:order}.} Recall that $\hat z_{k}= \ovx_k + \frac{\si}{\al}\ovw_k$, and define $\ep_k $ by $T_{[0,2]}\cro{\hat z_{k}}  =  \ovx_k + \ep_k$.
 Then $\E \cro{ T_{[0,2]}\cro{\hat z_{k}}}-p_k=\E \cro{ \ep_k}$ and it suffices to bound $|\E \cro{ \ep_k}|$.
Introducing the event $A= \{  |\frac{\si}{\al}\ovw_{k}| <  \ovx_k  \}$ and the complementary event $A^c$, we note first that $A \subseteq  \{\hat z_k \in [0,2]\}$ and thus $\ep_k =  \frac \si \al \hat w_k$ on $A$.
We have
\begin{align*}
  |\E \cro{ \ep_k}|
  \leq   |\E \cro{\ep_k \mathbbm{1}_A }|
  +  |\E \cro{\ep_k\mathbbm{1}_{A^c} }|
  &=  |\E \cro{\frac{\si}{\al}\ovw_{k}\mathbbm{1}_{A} }|
  +  |\E \cro{ \ep_k\mathbbm{1}_{A^c}}| \\
  &=  |\E \E \cro{\frac{\si}{\al}\ovw_{k}\mathbbm{1}_{A} \Big{|} \ovx_k }|
  +  |\E \cro{ \ep_k\mathbbm{1}_{A^c}}| \\
   &= |\E \cro{ \ep_k\mathbbm{1}_{A^c}}|
\end{align*}since $\ovw_k$ is a centered and symmetric random variable that is independent of $\ovx_k$.  Using the event $B=\{ 2 p_k \geq \ovx_k \geq p_k/2\}$ and the complementary event $B^c$, we have
\begin{align*}
  |\E \cro{ \ep_k\mathbbm{1}_{A^c}}|  \leq \E \cro{|\ep_k| \mathbbm{1}_{A^c \cap B} }| +  \E \cro{|\ep_k|\mathbbm{1}_{A^c\cap B^c} } &\leq   \E \cro{|\ep_k| \mathbbm{1}_{\{\frac \si \al|\ovw_{k}|\geq \frac 12 p_k \}} } +  2\E \cro{\mathbbm{1}_{B^c} } \\
  &\leq   \E \cro{\frac{\si}{\al} |\ovw_{k}| \mathbbm{1}_{\{\frac \si \al|\ovw_{k}|\geq \frac 12 p_k \}} } + 4e^{- \frac 18 n p_k }\\
  &=  2 p_k^{-1}\pa{ \E \cro{\frac{p_k}{2} |\frac{\si}{\al}\ovw_{k}| \mathbbm{1}_{\{\frac \si \al|\ovw_{k}|\geq \frac 12 p_k \}} } +  2p_k e^{- \frac 18 n p_k}}\\
   &\leq  2 p_k^{-1}\pa{ \E \cro{|\frac{\si}{\al}\ovw_{k}|^2  } +  2p_k e^{- \frac 18 n p_k }}
\end{align*}
where  we invoked (\ref{concentratoinBinomialJiao}-\ref{concentratoinBinomialJiao:bis}) from Lemma \ref{hatZ:bigP:nonInterac} in the second line. Then, by \eqref{eq2:lemmaC2} from Lemma \ref{LaplaceControl},
\begin{equation*}
  |\E \cro{ \ep_k\mathbbm{1}_{A^c}}| \leq 2 p_k^{-1} \pa{ \frac{\si^2 C_{L,2}}{\al^2 n} +  2p_k e^{-n p_k/8}}  \leq  2 p_k^{-1} \pa{ \frac{\si^2 C_{L,2}}{\al^2 n} +  \frac{16\ga}{en}}
\end{equation*}
where we used $x e^{-c n x} \leq \frac{\ga}{cen}$ for any $x\in [0,1]$ and any $c>0$. This concludes the proof of Lemma~\ref{lem:control:first:order}. \hfill $\square$

\medskip

\begin{lem}\label{hoeffding:zk}
For any $\al \in (0,\infty)$, and integers $K,n$ satisfying $ n \geq \log(Kn^{1/3})$, %\textcolor{red}{ajouter cette hypothèse sur $K,n$ au THM 3.4} 
we have
\begin{equation*}
      \Pb \left(|\hat{z}_k - p_k | > 96 \si \sqrt{\frac{ \log(Kn^{1/3})}{(\al^2\wedge 1)n}}\right)  \leq \frac{6}{K^3n}     \enspace.
\end{equation*}
\end{lem}

\medskip

\noindent \textbf{Proof of Lemma \ref{hoeffding:zk}.} Denoting $\de = c_1 \si \sqrt{\frac{ \log(Kn^{1/3})}{(\al^2\wedge 1) n}}$ with $c_1 \geq 1$ a numerical constant to be set later, we get from \eqref{hoeffding:binom} in Lemma \ref{hatZ:bigP:nonInterac} and \eqref{eq1:lemmaC2} in Lemma \ref{LaplaceControl} that
\begin{align*}
\Pb (|\hat{z}_k - p_k | > \de )  \leq \Pb(|\ovx_k-p_k| > \frac{\de}{2}) + \Pb(\frac{\si|\ovw_k|}{\al}>  \frac{ \de}{2 })
    &\leq  2 \left( e^{-\frac{n \de^2}{2}} + e^{-\frac{n(\al \de/\si)^2 }{32}} + e^{-\frac{n(\al \de/\si) }{8}  } \right)  \\
     &\leq  6\,  e^{-\frac{c_1  \log(Kn^{1/3}) }{32 }}
\end{align*}
which is upper bounded by $6/(K^3n)$ for $c_1= 96$. Lemma \ref{hoeffding:zk} is proved. \hfill $\square$

\bigskip
%%%%%%%%%%%%%%%%%%%%%%%%%%%%%%%%%%%%

\begin{lem}\label{lem:negassociattion}
We have \, $ \textup{Cov}\big{(}\hat{F}_{\ga}(k), \hat{F}_{\ga}(k') \big{)} \leq 0$ for any $k, k' \in[K]$, $k \neq k'$, and any $\ga >0$.
\end{lem}

\medskip

\textbf{Proof of Lemma \ref{lem:negassociattion}.}  We first state the definition of the negative association property.

\textit{Definition} (See \cite{negass}) Random variables $u_1,\ldots,u_K$ are said to be negatively associated (NA) if for every pair of disjoint subsets $A_1,A_2$ of $\{1,\ldots,K\}$, and any component-wise increasing functions $f_1,f_2$,
\begin{equation}\label{def:NA}
    \textup{Cov}\big{(}f_1(u_i, i\in A_1), f_2(u_j, j\in A_2)\big{)} \leq 0 \enspace.
\end{equation}

By corollary 5 of \cite{Jiao_2017}, random variables that are drawn from  a multinomial distribution, are NA. 
Hence, the  random variables $\hat{X}= (\hat{x}_1,\ldots,\hat{x}_K)$ are NA since $(\hat{x}_1,\ldots,\hat{x}_K)$ follows a multinomial distribution $\sim \mathcal{M}(n; (p_k)_{k\in[K]})$. Besides, the $\hat{W}= (\hat{w}_k)_{k\in [K]}$ are NA, as any set of independent random variables are NA \cite{negass}. Then, we get that  $(\hat{X}, \hat{W}) = (\hat{x}_1,\ldots,\hat{x}_K, \hat{w}_1,\ldots,\hat{w}_K)$ are NA since a standard closure property of NA is that the union of two independent sets of NA random variables is NA \cite{negass}. We can therefore use the definition \eqref{def:NA} of NA random variables to have 
\begin{equation*}
    \textup{Cov}\big{(}f_k(\hat{X},\hat{W}), f_{k'} (\hat{X},\hat{W})\big{)} \leq 0 \enspace, \qquad \forall k,k'\in [K],  k \neq k'
\end{equation*}
for  $f_k[(\hat{x}_1,\ldots,\hat{x}_K, \hat{w}_1,\ldots,\hat{w}_K)]$ $=\left[T_{[0,2]} (\hat{x}_k + \si \hat{w}_k /\al)  \right]^{\ga}$, which are component-wise increasing functions. The proof of Lemma \ref{lem:negassociattion} is complete. \hfill $\square$

%%%%%%%%%%%%%%%%%%%%%%%%%%%%%%%%%%%%%%%%%%%

%Sub-Section Small pk

%%%%%%%%%%%%%%%%%%
%%%%%%%%%%%%%%%%%%

\subsection{Bias and Variance on small values of $p_k$}
\begin{lem}\label{lem:bias:small:pk} Let $\ga, \al\in (0,\infty)$ and $k \in[K]$ and $c > 1$ be any numerical constant.  If $p_k \leq c / \sqrt{(\al^2 \wedge 1) n}$, then
\begin{equation*}
    \left| \E \hat{F}_{\ga}(k) - p_k^{\ga} \right| \leq  \frac{C}{((\al^2 \wedge 1) n)^{\ga/2}} \enspace,
\end{equation*}
\begin{equation*}
 \Va \pa{  \hat{F}_{\ga}(k) }   \leq   \frac{C' }{((\al^2 \wedge 1) n)^{\ga}} \enspace,
\end{equation*}
where $C, C'$ are constants depending only on $\ga$ and $c$.
\end{lem}

\medskip
%%%%%%%%%%%%%%%%%%%%%%%%%%%%%%%%%%%%%%%%

\noindent \textbf{Proof of Lemma \ref{lem:bias:small:pk}. } Recall that $\hat{F}_{\ga}(k) = \left(T_{[0,2]}\cro{\hat z_{k}} \right)^{\ga}$.
We have for any $s=1,2$,
\begin{equation*}
 \E \cro{ (\hat{F}_{\ga}(k))^s} =  \E \cro{ \left(T_{[0,2]}\cro{\hat z_{k}} \right)^{s\ga}} \leq  \E \cro{ |\hat z_{k}|^{s\ga}} \leq \frac{2^{s\ga} C_{BL,s\ga}}{((\al^2 \wedge 1)n)^{s\ga/2}} +  2^{s\ga} p_k^{s\ga}
\end{equation*}
using Lemma \ref{concentrationBinomialEtLaplace}. Then, we take $s=1$ to obtain the first bound announced in the lemma:
\begin{align*}
\left| \E \cro{\hat{F}_{\ga}(k)} - p_k^{\ga} \right| \leq   \E \cro{ \hat{F}_{\ga}(k)} +  p_k^{\ga} &\leq \frac{2^{\ga} C_{BL,\ga}}{((\al^2 \wedge 1)n)^{\ga/2}} +  (2^{\ga}+1) p_k^{\ga}\\
&\leq \frac{2^{\ga} C_{BL,\ga} + (2^{\ga}+1) c^{\ga}}{((\al^2 \wedge 1)n)^{\ga/2}}
\end{align*}
since $p_k \leq c / \sqrt{(\al^2 \wedge 1) n}$.
We finally take $s=2$ to get the second bound of the lemma:
\begin{align*}
\Va \pa{  \hat{F}_{\ga}(k) } \leq  \E \cro{ \hat{F}_{\ga}(k)^2 } \leq \frac{2^{2\ga} C_{BL,2\ga}+ 2^{2\ga} c^{2\ga}}{((\al^2 \wedge 1)n)^{\ga}}  \enspace .
\end{align*}
Lemma \ref{lem:bias:small:pk} is proved. \hfill $\square$

%%%%%%%%%%

%Sub Section Large pk

%%%%%%%%%%%%%%%%%%%%

\subsection{Bias and Variance on large values of $p_k$}  \label{section:append:largepK}

\begin{lem}\label{lem:bias:big:pk:first}For any $\ga, \al\in (0,\infty)$ and $k \in[K]$ with $p_k \in (0,1]$,  we have
\begin{equation*}
  \left|\E \cro{\hat{F}_{\ga}(k)^s} - p_k^{s\ga}\right| \leq  C \pa{ p_k^{s\ga}   e^{-\frac{n}{128 \si^2} \pa{(\al \wedge 1) p_k}^2 } + \frac{\mathbbm{1}_{\{s \ga \geq 2\}}}{((\al^2 \wedge 1) n)^{s\ga/2}}  + \frac{p_k^{s\ga -2}}{(\al^2 \wedge 1) n }  }, \enspace \forall s=1,2,
\end{equation*}where $C$ is a constant depending only on $\ga$.
\end{lem}

\medskip

%\noindent \textbf{Proof of Lemma \ref{lem:bias:big:pk:first}. }
The proof of Lemma \ref{lem:bias:big:pk:first} is inspired by the variance bound \cite[Lemma 28]{Jiao_2017} as it is based on Taylor's formula with the second derivatives of $x^{\ga}$ and $x^{2\ga}$. However, the result in \cite{Jiao_2017} holds for $\ga \in (0,1)$ in the case of direct observations (no privacy), whereas Lemma \ref{lem:bias:big:pk:first} holds for any $\ga >0$ in the case of sanitized observations (privacy). %The extension to any $\ga >0$ is possible as soon as the event $\hat{z}_k \in [p_k/2,  2 p_k]$ occurs with high probability, which is guaranteed by the assumption of Lemma \ref{lem:bias:big:pk:first} (i.e. $p_k > c / \sqrt{(\al^2 \wedge 1) n}$ for some large enough numerical constant $c$).
We postpone the (relatively long) proof to the end of section \ref{section:append:largepK}.

\bigskip

\begin{lem}\label{lem:bias:big:pk} Let $\ga, \al\in (0,\infty)$ and $k \in[K]$ and $c>0$ be any numerical constant. If $p_k \geq c/ \sqrt{(\al^2 \wedge 1) n}$, then
\begin{equation*}
  \left|\E \cro{\hat{F}_{\ga}(k)^s} - p_k^{s\ga}\right| \leq  C \pa{  \frac{1}{((\al^2 \wedge 1) n)^{s\ga/2}}  + \mathbbm{1}_{\{ s \ga \geq 2\}}\frac{p_k^{s\ga -2}}{(\al^2 \wedge 1) n }  }, \enspace \forall s=1,2,
\end{equation*}where $C$ is a constant depending only on $\ga$ and $c$.
\end{lem}

\medskip

\noindent \textbf{Proof of Lemma \ref{lem:bias:big:pk}.} We invoke Lemma \ref{lem:bias:big:pk:first}. We bound the first error term
\begin{equation*}
    p_k^{s\ga}   e^{-\frac{n}{128 \si^2} \pa{(\al \wedge 1) p_k}^2 } \leq \pa{\frac{64  \si^2 s\ga}{(\al^2 \wedge 1) e n}}^{s\ga/2}
\end{equation*}
where we used  $ x^{s\ga} e^{-c n x^2} \leq \pa{\frac{s\ga}{2cen}}^{s\ga/2}$ for $x \in [0,1]$ and any $c>0$.
%The second error term  of Lemma \ref{lem:bias:big:pk:first} is smaller than $1/ ((\al^2 \wedge 1) n)^{s\ga/2}$.
The third error term of Lemma \ref{lem:bias:big:pk:first} satisfies, for $s\ga \in (0,2)$
\begin{equation}\label{lemVarrepeatNo}
   \mathbbm{1}_{\{ s \ga \in (0,2)\}} \frac{p_k^{s\ga -2}}{(\al^2 \wedge 1) n } \leq   \frac{\Big{(}\frac{c}{\sqrt{(\al^2 \wedge 1)n}}\Big{)}^{s\ga -2} }{(\al^2 \wedge 1) n } \leq \frac{c^{s\ga -2} }{((\al^2 \wedge 1) n)^{s \ga/2} }
\end{equation} since $p_k \geq  c/ \sqrt{(\al^2 \wedge 1)n}$. The proof of Lemma \ref{lem:bias:big:pk} is complete. \hfill $\square$

\bigskip

%%%%%%

\begin{lem}\label{lem:var:big:pk} Under the assumptions of Lemma \ref{lem:bias:big:pk}, we have
\begin{equation*}
       \Va \pa{  \hat{F}_{\ga}(k) } \leq   C \bigg{(} \frac{1}{((\al^2 \wedge 1) n)^{\ga}}  + \mathbbm{1}_{\{  \ga \geq 1\}}\frac{p_k^{2\ga -2}}{(\al^2 \wedge 1) n } \bigg{)}
\end{equation*}
for a constant $C$ depending only on $\ga$ (and $c$).
\end{lem}

\medskip

\noindent \textbf{Proof of Lemma \ref{lem:var:big:pk}. } We have, similarly to \cite{Jiao_2017},
\begin{align}\label{Var:nonInter:proof:BigP}
    \Va \pa{  \hat{F}_{\ga}(k) } &= \E \cro{\hat{F}_{\ga}(k)^2} - \pa{\E \hat{F}_{\ga}(k)}^2 = \E \cro{\hat{F}_{\ga}(k)^2} - p_k^{2\ga} + p_k^{2\ga} - \pa{\E \hat{F}_{\ga}(k)}^2 \nonumber \\
    &\leq \left|\E \cro{\hat{F}_{\ga}(k)^2} - p_k^{2\ga}\right| + \left|p_k^{2\ga} - \pa{\E \hat{F}_{\ga}(k) - p_k^{\ga} + p_k^{\ga}}^2\right| \nonumber \\
    &\leq  \left|\E \cro{\hat{F}_{\ga}(k)^2} - p_k^{2\ga}\right| + \left|\E \hat{F}_{\ga}(k) - p_k^{\ga}\right|^2 + 2p_k^{\ga} \left|\E \hat{F}_{\ga}(k) - p_k^{\ga}\right| \enspace.
\end{align}
Using Lemma \ref{lem:bias:big:pk} to bound the two first terms of \eqref{Var:nonInter:proof:BigP}, and Lemma \ref{lem:bias:big:pk:first} for the last term, we get
\begin{align}\label{6terms}
     \Va \pa{  \hat{F}_{\ga}(k) } \leq   C \bigg{(} &\frac{1}{((\al^2 \wedge 1) n)^{\ga}}  + \mathbbm{1}_{\{  \ga \geq 1\}}\frac{p_k^{2\ga -2}}{(\al^2 \wedge 1) n }\nonumber \\
      &+   \frac{1}{((\al^2 \wedge 1) n)^{\ga}}  + \mathbbm{1}_{\{ \ga \geq 2\}}\frac{ p_k^{2(\ga -2)}}{((\al^2 \wedge 1) n)^2 }  \\
       &+ 2p_k^{2\ga}   e^{-\frac{n}{128 \si^2} \pa{(\al \wedge 1) p_k}^2 } + \frac{2p_k^{\ga}\mathbbm{1}_{\{ \ga \geq 2\}}}{((\al^2 \wedge 1) n)^{\ga/2}}  + \frac{2 p_k^{2\ga -2}}{(\al^2 \wedge 1) n }  \bigg{)} \nonumber .
\end{align} We bound the fifth term of \eqref{6terms}:
\begin{equation*}
    2p_k^{2\ga}   e^{-\frac{n}{128 \si^2} \pa{(\al \wedge 1) p_k}^2 } \leq 2 \pa{\frac{128 \si^2 \ga}{(\al^2\wedge 1)en}}^{\ga}
\end{equation*}using  $ x^{2\ga} e^{-c' n x^2} \leq \pa{\frac{\ga}{c'en}}^{\ga}$ for any $x \in [0,1]$ and any $c'>0$. Hence, the first, third and fifth terms of \eqref{6terms} are  of the order of $((\al^2\wedge 1)n)^{-\ga}$ at most.
We now bound the fourth term of \eqref{6terms} using  $p_k\geq c / \sqrt{\al^2 \wedge 1) n}$ :
\begin{equation*}
    \frac{ p_k^{2(\ga -2)}}{((\al^2 \wedge 1) n)^2 } =  \frac{ p_k^{2\ga -2} p_k^{-2}}{((\al^2 \wedge 1) n)^2 } \leq  \frac{ p_k^{2\ga -2}}{c^2(\al^2 \wedge 1) n }
\end{equation*} and similarly the sixth term of \eqref{6terms}:
\begin{equation*}
    \frac{2p_k^{\ga}\mathbbm{1}_{\{ \ga \geq 2\}}}{((\al^2 \wedge 1) n)^{1 +(\ga/2) -1}} \leq \frac{2p_k^{\ga} (p_k/c)^{\ga-2}\mathbbm{1}_{\{ \ga \geq 2\}}}{(\al^2 \wedge 1) n} =  \frac{ 2p_k^{2\ga -2} \mathbbm{1}_{\{ \ga \geq 2\}}}{c^{\ga-2}(\al^2 \wedge 1) n } \enspace.
\end{equation*} %Hence, the second, fourth and sixth terms of \eqref{6terms} are  upper bounded by  $\mathbbm{1}_{\{ \ga \geq 1 \}}p_k^{2\ga -2}/(\al^2\wedge 1)n $. 
Hence, we have the desired bound for the second, fourth and sixth terms of \eqref{6terms}. Finally, for the last term of \eqref{6terms} we have
\begin{equation*}
    \frac{ p_k^{2\ga -2}}{(\al^2 \wedge 1) n }=  \frac{ p_k^{2\ga -2}\mathbbm{1}_{\{  \ga \in (0,1)\}}}{(\al^2 \wedge 1) n } +    \frac{ p_k^{2\ga -2}\mathbbm{1}_{\{  \ga \geq 1\}}}{(\al^2 \wedge 1) n } \leq \frac{2c^{2\ga -2} }{((\al^2 \wedge 1) n)^{ \ga} } +   \frac{ p_k^{2\ga -2}\mathbbm{1}_{\{  \ga \geq 1\}}}{(\al^2 \wedge 1) n }
\end{equation*}
using  \eqref{lemVarrepeatNo} for $s=2$.
This concludes the proof of of Lemma \ref{lem:var:big:pk}. \hfill $\square$

\bigskip
%%%%%%%%%%%%%%%%%%%%%%%%%%

%STOOOOOOOOOOOOOOOOOOOOOOOOOOOOP

\noindent \textbf{Proof of Lemma \ref{lem:bias:big:pk:first}. }
Denoting $f_s(x)=x^{s\ga}$ for $s=1,2$, and $Y= T_{[0,2]}\cro{\hat z_{k}}$, we have by Taylor's formula,
\begin{equation}\label{resteintegral}
    f_s(Y) = f_s(p_k) + f'_s(p_k)(Y-p_k)+  R(Y, p_k)
    \end{equation}
where the remainder is defined by
    \begin{equation}\label{resteintegral:bis}
    R(Y, p_k) = \int_{p_k}^Y (Y-w)f_s''(w)dw = \frac{1}{2}f_s''(w_Y) (Y-p_k)^2
    \end{equation}
where $w_Y$ lies between $Y$ and $p_k$. We get
\begin{equation}\label{biasEqualReste}
    |\E f_s(Y) - f_s(p_k)| \leq | \E R(Y, p_k)| + |\E f_s'(p_k)(Y-p_k) | \enspace .
\end{equation}
Thus, to prove the lemma, it suffices to bound the remainder $|\E R(Y, p_k)|$ and  the first order term $ |\E f_s'(p_k)(Y-p_k) |$. We control the latter using Lemma \ref{lem:control:first:order},
\begin{align*}
    |\E f_s'(p_k)(Y-p_k) | = s\ga p_k^{s\ga-1}|\E (Y-p_k) | \leq \frac{2s\ga p^{s\ga -2}}{(\al^2 \wedge 1)n} \pa{ \si^2 C_{L,2} +  \frac{16\ga}{e}} \enspace.
\end{align*} For the remainder, we use the decomposition
\begin{equation}\label{eq1:proof:bias:DL2}
    | \E R(Y, p_k)| \leq    \E \cro{|R(Y, p_k)| \mathbbm{1}(Y< p_k/2)} + \E \cro{|R(Y, p_k)| \mathbbm{1}(Y\geq p_k/2)}
\end{equation}
and we bound separately the two terms of the RHS.

\textit{$1^{\circ}.$ First term in the RHS of \eqref{eq1:proof:bias:DL2}. }
\begin{align*}
\E \cro{|R(Y, p_k)| \mathbbm{1}(Y< p_k/2)}  &\leq    \underset{y\leq p_k/2}{\textup{sup}}|R(y, p_k)|  \E \cro{\mathbbm{1}(Y< p_k/2) } \\
&=   \underset{y\leq p_k/2}{\textup{sup}}|R(y, p_k)|  \E \cro{\mathbbm{1}( \hat{z}_k< p_k/2)  } \\
     &\leq  \underset{y\leq p_k/2}{\textup{sup}}|R(y, p_k)|\, 3\,  e^{-\frac{n}{128 \si^2} \pa{(\al \wedge 1) p_k}^2 }
\end{align*}using Lemma \ref{concentrationBinomialEtLaplace}.
We control $R(y, p_k)$  for any $y \in [0,p_k/2]$,
\begin{align*}
   |R (y, p_k)| & \leq \int_y^{p_k} (w-y) |f_s''(w)| dw \leq \int_y^{p_k} (w-y)s\ga |s\ga -1| w^{s\ga-2}dw\\
   &\leq s\ga |s\ga -1| \int_y^{p_k} w^{s\ga-1}dw \leq s\ga |s\ga -1| \int_0^{p_k} w^{s\ga-1}dw = |s\ga -1|  p_k^{s\ga}\enspace.
\end{align*}
We gather the last two displays to get
\begin{equation*}
   \E \cro{|R(Y, p_k)| \mathbbm{1}(Y< p_k/2)} \leq  |s\ga -1| p_k^{s\ga} \, 3\,  e^{-\frac{n}{128 \si^2} \pa{(\al \wedge 1) p_k}^2 }  \enspace.
\end{equation*}

\textit{$2^{\circ}$. Second term in the RHS of \eqref{eq1:proof:bias:DL2}. } We separate our analysis in two different ranges of values of $\ga$.

\textit{$2^{\circ}.1$. Case $s\ga \in (0,2)$: } Starting from \eqref{resteintegral:bis} we have
\begin{align}\label{nonrepeat}
    \E \cro{|R(Y, p_k)| \mathbbm{1}(Y\geq p_k/2)}& = \frac{ s\ga |s\ga -1|}{2}\E \cro{w_Y^{s\ga -2} (Y-p_k)^2 \mathbbm{1}(Y\geq p_k/2)}  \\
    &\leq \frac{ s\ga |s\ga -1|}{2}  \big{(}\frac{p_k}{2}\big{)}^{s\ga-2} \E \cro{(Y-p_k)^2 }  \nonumber \\
    &\leq  s\ga |s\ga -1|  2^{1-s\ga } p_k^{s\ga -2} \frac{ C_{BL,2}}{(\al^2 \wedge 1) n } \nonumber
\end{align}
where   we used $\E \cro{(Y-p_k)^2 } \leq \E\cro{   (\hat{z}_k-p_k)^2  }$ and Lemma \ref{concentrationBinomialEtLaplace}.

\smallskip

\textit{$2^{\circ}.2$. Case $s\ga \geq 2$: } A plug of $w_Y^{s\ga-2} \leq p_k^{s\ga-2} + Y^{s\ga-2}$ into \eqref{nonrepeat} gives
\begin{equation}\label{ptiteEqNew}
    \E \cro{|R(Y, p_k)| \mathbbm{1}(Y\geq p_k/2)} \leq  \frac{ s\ga |s\ga -1|}{2}\E \cro{(p_k^{s\ga-2}+Y^{s\ga-2}) (Y-p_k)^2 \mathbbm{1}( Y\geq p_k/2) } .
\end{equation}
We bound the first part of \eqref{ptiteEqNew} as in \eqref{nonrepeat},
\begin{equation*}
    \E \cro{p_k^{s\ga-2}(Y-p_k)^2 \mathbbm{1}( Y\geq p_k/2) }   \leq p_k^{s\ga-2} \frac{ C_{BL, 2}}{((\al^2 \wedge 1)n)}  \enspace.
\end{equation*}
For the second part of \eqref{ptiteEqNew}, we get from Cauchy-Schwarz  that
\begin{align*}
\E \cro{Y^{s\ga -2} (Y-p_k)^2 \mathbbm{1}( Y\geq 2p_k)}
      &\leq \E \cro{Y^{2(s\ga -2)}}^{1/2} \E \cro{ (Y-p_k)^4 }^{1/2} \\
      &\leq \pa{\frac{2^{2(s\ga -2)} C_{BL,2(s\ga -2)}}{((\al^2 \wedge 1)n)^{s\ga -2}} +  2^{2(s\ga -2)} p_k^{2(s\ga -2)}}^{1/2} \pa{\frac{C_{BL,4}}{((\al^2 \wedge 1)n)^{2}}}^{1/2}\\
      &\leq \pa{\frac{2^{s\ga -2} \sqrt{C_{BL,2(s\ga -2)}}}{((\al^2 \wedge 1)n)^{(s\ga -2)/2}} +  2^{s\ga -2} p_k^{s\ga -2}} \frac{\sqrt{C_{BL,4}}}{(\al^2 \wedge 1)n}
\end{align*}
where in the second inequality we used $\E \cro{ Y^{2r} } \leq \E \cro{ \hat{z}_k^{2r}}$ and $\E \cro{ (Y-p_k)^{2r} } \leq \E \cro{ (\hat{z}_k-p_k)^{2r}}$ for any $r>0$ and Lemma \ref{concentrationBinomialEtLaplace}; in the third inequality we used $\sqrt{a+b} \leq \sqrt{a} + \sqrt{b}$ for any $a,b >0$. A plug of the last two displays into \eqref{ptiteEqNew} concludes the case $s\ga \geq 2$. 

Going back to \eqref{eq1:proof:bias:DL2}, we have bounded the remainder $\E R(Y, p_k)$. Lemma \ref{lem:bias:big:pk:first} is proved. \hfill $\square$

%%%%%%%%%%%%%%%%%%%%%.

%SECTION Auxiliary Lemmas

%%%%%%%%%%%%%%%%%%%%%%%%%%%%

 \section{Auxiliary lemmas for upper bounds }

\begin{lem}\label{hatZ:bigP:nonInterac}
Let $p\in(0,1]$, and  $x_1,\ldots,x_n \overset{iid}{\sim} \textup{B}(p)$ be independent Bernoulli random variables with parameter $p$. Then, the mean $\ovx = \frac{1}{n} \sum_{i=1}^n x_i$  satisfies, for any $\de >0$ ,
\begin{equation}\label{concentratoinBinomialJiao}
\P \pa{   \ovx \leq (1-\de)p } \leq e^{- \frac{ \de^2 n p}{ 2}} \enspace,
\end{equation}
\begin{equation}\label{concentratoinBinomialJiao:bis}
\P \pa{   \ovx \geq (1+\de)p } \leq e^{- \frac{ \de^2 n p}{ 2 + \de}} \enspace,
\end{equation}
and
\begin{equation}\label{hoeffding:binom}
     \Pb (|\hat x -p| \geq \de)  \leq 2 e^{-2 \de^2 n}  \enspace .
\end{equation}
We also have, for any $r>0$,
\begin{equation}\label{ep:Mean:Binom}
    \E\cro{ |\ovx -p|^{r }} \leq \frac{C_{B,r}}{n^{r/2}}
\end{equation}
where $C_{B, r}$ is a constant depending only on $r$.
\end{lem}

\medskip

\noindent \textbf{Proof of Lemma \ref{hatZ:bigP:nonInterac}.}  The concentration inequalities (\ref{concentratoinBinomialJiao}-\ref{concentratoinBinomialJiao:bis}) are  one form of Chernoff bounds. The control \eqref{hoeffding:binom} is Hoeffding's inequality applied to i.i.d Bernoulli random variables. Finally, for (\ref{ep:Mean:Binom}), see \cite{Wainwright} or adapt the proof of Lemma \ref{LaplaceControl} below. \hfill $\square$

\bigskip

\begin{lem}\label{LaplaceControl}
Let $w_1,\ldots,w_n \overset{iid}{\sim} \textup{L}(1)$ be independent Laplace random variables with parameter $1$. Denoting the mean by $\ovw = \frac{1}{n} \sum_{i=1}^n w_i$ , we have
\begin{align}\label{eq1:lemmaC2}
   \Pb (\ovw >t) \vee \Pb (\ovw < -t) &\leq \exp \cro{-\frac{n}{2}(\frac{t^2}{4} \wedge \frac{t}{2}) } \nonumber \\
   &\leq \exp \cro{-\frac{n}{8}t^2 }  + \exp \cro{-\frac{n}{4} t } \enspace.
\end{align}
Besides, for any real $r>0$, there exists a constant $C_{L,r} \geq  1$, depending only on $r$, such that
\begin{equation}\label{eq2:lemmaC2}
    \E \pa{  \left| \ovw \right|^{r}} \leq \frac{C_{L,r}}{n^{r/2}} \enspace.
\end{equation}
\end{lem}

\medskip

\noindent \textbf{Proof of Lemma \ref{LaplaceControl}.} A random variable $x$ is said to be sub-exponential with parameter $\lambda$, denoted $x \sim \textup{subE}(\lambda)$, if $\E x = 0$ and its moment generating function satisfies
\begin{equation*}
    \E [ e^{s x} ] \leq e^{\lambda^2 s^2 /2}, \quad  \forall  \, |s| < \frac{1}{\lambda}.
\end{equation*}

Let $x_1,\ldots,x_n$ be independent random variables such that  $x_i \sim \textup{subE}(\lambda)$.  Bernstein's inequality \cite{Wainwright} entails that, for any $t>0$, the mean $\hat x = \frac{1}{n} \sum_{i=1}^n x_i$ satisfies
\begin{equation}\label{proof:bernstein:exp}
    \Pb (\hat x >t) \vee \Pb (\hat x < -t) \leq \exp \cro{-\frac{n}{2}(\frac{t^2}{\lambda^2} \wedge \frac{t}{\lambda}) } \enspace.
\end{equation}
Then, for any real $r >0$ we have
\begin{align*}
    \E |\ovx | = \int_0^{\infty} \Pb(|\ovx |^{r} > t) dt =  \int_0^{\infty} \Pb(|\ovx | > t^{1/r}) dt  \leq  \int_0^{\infty} 2 e^{ -\frac{n t^{2/r}}{2 \lambda^2}}  dt + \int_0^{\infty} 2 e^{ -\frac{n t^{1/r}}{2 \lambda}}  dt
\end{align*} so that, using $u= \frac{n t^{2/r}}{2 \lambda^2}$ and $v =\frac{n t^{1/r}}{2 \lambda}$,
\begin{align}\label{proof:mean!subExp}
    \E |\ovx | & \leq   \pa{\frac{2 \lambda ^2}{n}}^{r/2} r \int_0^{\infty}  e^{ -u} u^{(r /2) -1}  du \ \, + \ \,  2 \pa{\frac{2 \lambda }{n}}^{r} r \int_0^{\infty}  e^{ -v} v^{r  -1} dv \nonumber  \\
    &= \pa{\frac{2 \lambda ^2}{n}}^{r/2} r \Gamma(r/2) \ \, + \ \,  2 \pa{\frac{2 \lambda }{n}}^{r} r \Gamma(r) \nonumber \\
    &\leq 2^{r +2} \lambda^{r}  r \left[ \Gamma(r/2) + \Gamma(r) \right] \frac{1}{n^{r/2}} \enspace.
\end{align}

Let $w \sim \textup{L}(1)$ be a random variable of Laplace distribution with parameter $1$. Observe that $\Pb (|w| > t ) = e^{-t}$ for $t\geq 0$, and
\begin{equation*}\E [ e^{s w} ] \leq e^{2 s^2}, \quad \textup{ if } \, |s| < \frac{1}{2} .\end{equation*}
Hence,
$w$ is sub-exponential with parameter $2$, i.e. $w \sim \textup{subE}(2)$. We can take $\lambda = 2$ in (\ref{proof:bernstein:exp}-\ref{proof:mean!subExp}) to conclude the proof of Lemma \ref{LaplaceControl}, choosing $C_{L,r} = 2^{2r +2}  r \left[ \Gamma(r/2) + \Gamma(r) \right]$. \hfill $\square$

%%%%%%%%%%%%%%%%%%

%%%%%%%%%%%%%%%%%%

%New Section

%%%%%%%%%%%%%%%%

\section{Proofs of lower bounds}
\label{section:proof:LB}

\noindent \textbf{Proof of Proposition 2.2. } Recall that $\hat{z}_k = \frac{1}{n}\sum_{i=1}^n z_{ik}$, where $z_{ik} = \mathbbm{1}_{\{x_i = k\}} +\frac{\si}{\alpha} \cdot w_{ik}$, with $\E z_{ik} = p_k$ and $\textup{Var}(z_{ik}) = p_k(1-p_k) + \frac{2\si^2}{\al^2}$. Note that $\tilde \tau :=  \frac{\si}{\sqrt{\al^2 n}}$ lies in $[0,2]$, and that $\textup{Var}(z_{ik}) \geq (\sqrt{n}  \tilde \tau)^2$. By the central limit theorem, $ \sqrt{n} \frac{\hat{z}_k-p_k}{\sqrt{\textup{Var}(z_{ik}) }}$ has an asymptotic  standard normal distribution, so we have $\Pb(\sqrt{n} \frac{\hat{z}_k-p_k}{\sqrt{\textup{Var}(z_{ik}) }}\geq 1) \geq c_1$ for some numerical  constant $c_1>0$ and $n$ large enough. We write $\hat{z}_k = \sqrt{n} \frac{\hat{z}_k-p_k}{\sqrt{\textup{Var}(z_{ik}) }} \cdot \frac{\sqrt{\textup{Var}(z_{ik})}}{\sqrt{n}} + p_k \geq \frac{\sqrt{\textup{Var}(z_{ik})}}{\sqrt{n}}$ with probability larger than $c_1$, thus leading to 
\begin{equation*}
    \E \cro{\left(T_{[0,2]}(\hat{z}_k)\right)^{\ga}} - p_k^{\ga} \geq c_1  \pa{T_{[0,2]}\Big{(}\frac{\sqrt{\textup{Var}(z_{ik})}}{\sqrt{n}}\Big{)}}^{\ga} - p_k^{\ga} = c_1 \tilde \tau^{\ga}  - p_k^{\ga} \geq \frac{c_1 \tilde \tau^{\ga}}{2} \enspace, \qquad \textup{ as } \ \, n \rightarrow \infty
\end{equation*}
for all $p_k \leq \pa{\frac{c_1 }{2}}^{1/\ga} \tilde \tau$. Denoting by $\Kcal_{\leq (c_1 /2)^{1/\ga} \tilde \tau}$ the number of such $p_k$ satisfying the latter inequality, we get 
\begin{equation}\label{proof:asymp:small:pk}
  \sum_{k\in \Kcal_{\leq (c_1 /2)^{1/\ga} \tilde \tau} }  \E \cro{ \left(T_{[0,2]}(\hat{z}_k)\right)^{\ga}} - p_k^{\ga}  \geq \frac{c_1 \tilde \tau^{\ga} |\Kcal_{\leq (c_1 /2)^{1/\ga} \tilde \tau}|}{2} \enspace, \qquad \textup{ as } \ \, n \rightarrow \infty \enspace.
\end{equation}
Hence, the lower bound announced in Proposition 2.2 holds in particular for any  $p=(p_1,\ldots,p_K) \in \Pcal_K$ such that $|\Kcal_{\leq (c_1 /2)^{1/\ga} \tilde \tau}| = K$. However, this last equality entails that  $K$ satisfies the following restriction $K\gtrsim_{\ga} (\tilde \tau)^{-1} \gtrsim_{\ga} \sqrt{\al^2 n}$ since $\sum_{k=0}^K p_k =1$. We remove this restriction in the sequel.

Let $C>0$ be some constant that will be set later, and that only depends on $\ga$. If $K \leq  C \left( 1 \vee (\al^2 n)^{\frac{\ga}{2}-\frac{1}{2}} \right)$, then the lower bound of  Proposition 2.2 follows directly from Theorem 2.6.   
We can therefore assume that
\begin{equation}\label{bornein:asymp!2volet}
    K\geq C \left( 1 \vee (\al^2 n)^{\frac{\ga}{2}-\frac{1}{2}} \right) \enspace.
\end{equation} 
Let $p=(p_1,\ldots,p_K) \in \mathcal{P}_K$ such that $p_j \leq  \pa{\frac{c_1 }{2}}^{1/\ga} \tilde \tau $ for all $j\in[K-1]$ , and $p_K \in[0,1]$ so that $\sum_{k=1}^K p_k =1$. By Lemma \ref{lem:bias:small:pk} and \ref{lem:bias:big:pk}, the bias of estimation of $p_K$ is bounded  by
\begin{equation*}
    \left|  \E \cro{\left(T_{[0,2]}(\hat{z}_K)\right)^{\ga}}   - p_K^{\ga} \right| \leq   C' \pa{  \frac{1}{(\al^2  n)^{\ga/2}}  + \mathbbm{1}_{\{ \ga \geq 2\}}\frac{1}{\al^2  n }  } \enspace,
\end{equation*}where $C'$ is a constant depending only on $\ga$. Combining with \eqref{proof:asymp:small:pk}, we get 
\begin{align*}
  \sum_{k=1 }^K  \E \left(T_{[0,2]}(\hat{z}_k)\right)^{\ga} - p_k^{\ga} & \geq \frac{c_1 \tilde \tau^{\ga} (K-1)}{2} -  \frac{C'}{(\al^2  n)^{\ga/2}} -  \mathbbm{1}_{\{ \ga \geq 2\}}\frac{C'}{\al^2  n }  \\
  &\geq \frac{c_1  K}{4 (\al^2  n)^{\ga/2}} -  \frac{C'}{(\al^2  n)^{\ga/2}} -  \mathbbm{1}_{\{ \ga \geq 2\}}\frac{C'}{\al^2  n } \enspace.
\end{align*} Hence, it suffices to choose a large enough constant $C$   in \eqref{bornein:asymp!2volet} to have 
\begin{equation*}
     \sum_{k=1 }^K  \E \left(T_{[0,2]}(\hat{z}_k)\right)^{\ga} - p_k^{\ga}\geq \frac{C''  K}{ (\al^2  n)^{\ga/2}}
\end{equation*}for some constant $C''$ depending only on $\ga$. We have proved the desired lower bound under the assumption \eqref{bornein:asymp!2volet}. 
The proof of Proposition 2.2 is complete. \hfill $\square$

\bigskip

\noindent \textbf{Proof of Theorem 2.6.}  Fix $\ga >0 , \ga \neq 1$. Let $\tilde \tau:= \frac{\tilde C}{\sqrt{\al^2 n}}$ for a constant $\tilde C \in (0,1)$ that will be set later, and which only depends on $\ga$. Let us start with the case  $K=2$. Define two probability vectors  $p = (p_1,p_2) = (1- \tilde \tau, \tilde \tau)$ and $q= (q_1,q_2) = (1-\tilde \tau/2, \tilde \tau/2)$. Then for a small enough constant $\tilde C$, we have
\begin{align*}
    \Delta := |F_{\ga}(p) - F_{\ga}(q)|& = \left|(1- \tilde \tau)^{\ga} - (1- \tilde \tau/2)^{\ga} + \tilde \tau^{\ga} -( \tilde \tau/2)^{\ga}\right| \\
    &= \left|-\frac{\ga \tilde \tau}{2} + O(\tilde \tau^2) + \tilde \tau^{\ga} (1-\frac{1}{2^{\ga}})\right|
\end{align*}
where we used $(1-x)^{\ga} = 1- \ga x + O(x^2)$ for any real $x\in(0,\tilde C)$. If $\ga\in(0,1)$, we can choose $\tilde C$ small enough to have 
\begin{align*}
    \Delta &= \tilde \tau^{\ga} \left|-\frac{\ga \tilde \tau^{1-\ga}}{2} + O(\tilde \tau^{2-\ga}) + (1-\frac{1}{2^{\ga}})\right| \geq C \tilde \tau^{\ga}
\end{align*}for some constant $C$ depending only on $\ga$. Similarly, if $\ga\ > 1$, we have
\begin{align*}
    \Delta &= \tilde \tau \left|-\frac{\ga }{2} + O(\tilde \tau) + \tilde \tau^{\ga-1}(1-\frac{1}{2^{\ga}})\right| \geq C \tilde \tau \enspace.
\end{align*}

For any $\al$-LDP mechanism $Q$, denote by $Qp$ and $Qq$ the measures corresponding to the channel $Q$ applied to the probability vectors $p$ and $q$. Corollary 3 of \cite{duchi2018minimax} ensures that the Kullback-Leibler divergence between $Qp$ and $Qq$ is bounded by 
\begin{equation*}
D_{kl}(Qp,Qq) \leq 4(e^{\al}-1)^2 n \left(d_{TV}(p,q)\right)^2\enspace,
\end{equation*}
i.e. by $n$ times the square of the total variation distance between $p$ and $q$, up to a constant depending on $\alpha$.
Then we have
\begin{equation}\label{proof:KL:lb:thm1}
D_{kl}(Qp,Qq) \leq  4(e^{\al}-1)^2 n  \left(\sum_{k=1}^2 |p_k-q_k|\right)^2 \leq  4(e^{\al}-1)^2 n \tilde \tau^2 \leq 36 \tilde C^2 
\end{equation}
where the last inequality follows from $e^x -1 \leq 3x$ for any $x\in [0,1]$.

For any vector $\theta =(\theta_1,\theta_2)$, $\theta_i \geq 0$, we denote the functional at $\theta$ by $F_{\ga}(\theta) =\sum_{k=1}^2 \theta_k^{\ga}$. We use a standard lower bound  method based on two hypotheses, see e.g. Theorem 2.1 and 2.2 in \cite{Tsybakov2009IntroductionTN}, to get for any  estimator $\hat{F}$,
\begin{equation*} \underset{\theta  \in \{p, q \}}{\textup{sup}} \, \P_{\theta}\left( |\hat{F} - F_{\ga}(\theta)| \geq \frac{\Delta}{2}\right) \geq \frac{1- \sqrt{D_{kl}(Qp,Qq)/2}}{2} \enspace.
\end{equation*}
Then we deduce from \eqref{proof:KL:lb:thm1} that
\begin{equation*} \underset{\theta  \in \{p, q \}}{\textup{sup}} \, \P_{\theta}\left( |\hat{F} - F_{\ga}(\theta)| \geq \frac{\Delta}{2}\right)  \geq \frac{1- 3\sqrt{2}\tilde C}{2} \geq \frac{1}{4} \enspace,
\end{equation*}
choosing $\tilde C \leq 1/(6\sqrt{2})$. We have proved the desired lower bound in the case $K=2$.

We can actually prove the same lower bound for any integer $K\geq 2$, with the following slight modification in the proof written above. Choose $p_k$, $q_k$, $k\geq 3$ such that $p_k=q_k$ and $p_k \leq \tilde{C}/(4K n)$. Then change the $p_1$ and $q_1$ above accordingly (to have probability vectors). This affects neither the order of the separation $\Delta$, nor the bound on the KL-divergence between the measures $Qp$ and $Qq$. This concludes the proof of Theorem 2.6.
\hfill $\square$
%%%%%%%%%%%%%
%%%%%%%%%%%%%%%%%%%%%%%%%%

\bigskip

\noindent \textbf{Proof of Theorem 2.7.} If $K < 4$, then the lower bounds are a direct consequence of Theorem 2.6. We assume therefore that $K\geq 4$. For the ease of exposition, we also assume that $K$ is even (the case of an odd $K$ being similar). Let $\tilde K$ be a   positive even  integer in $[K]$. Let $p=(p_1,\ldots,p_K)$ be any probability vector such that two consecutive coordinates are equal $p_{2k-1} = p_{2k}$ for  $k\in [\tilde K/2]$, and the remaining coordinates satisfy $p_k=p_{k'}$ for all $k,k'\geq \tilde K + 1$. Similarly, let $\de=(\de_1,\ldots,\de_{ K})$ be a vector of perturbations such that,  two consecutive perturbations are equal $\de_{2k-1} =  \de_{2k}$, $k\in [\tilde K/2]$,   and the others are equal to zero:  $\de_k =0$ , $\forall$ $k\geq \tilde K +1$. Each perturbation is smaller than (half of) the corresponding probability: $0 \leq \de_{k} \leq p_{k}/2$, $k\in[\tilde K]$. Given any $k\in[K/2]$ and any vector $q=(q_1,\ldots,q_{K})$, define the operator $T_k(q) = (0,\ldots,0,q_{2k-1},-q_{2k},0,\ldots,0)$. We are now ready to introduce the following collection of vectors $p^{(\nu)}$, $\nu \in \Vcal \{-1,1\}^{\tilde K/2}$:
\begin{align*}
    p^{(\nu)} &= p + \sum_{k=1}^{\tilde K/2} \nu_k T_k(\de) \\
    &= (p_1, p_2, p_3,p_4,  \ldots, p_{K-1}, p_{K}) + ( \nu_1\de_1,-\nu_1 \de_2, \ldots ,\nu_{\tilde K/2}\de_{\tilde K-1},-\nu_{\tilde K/2} \de_{\tilde K},0,\ldots,0) \\
    &= (p_2, p_2, p_4, p_4, \ldots, p_{\tilde K}, p_{ \tilde K}, p_K,\ldots,p_K) + ( \nu_1\de_2,-\nu_1 \de_2, \ldots ,\nu_{\tilde K/2}\de_{\tilde K},-\nu_{\tilde K/2} \de_{\tilde K},0,\ldots,0)  \enspace.
\end{align*}
Observe that each $p^{(\nu)}$, $\nu \in \Vcal \{-1,1\}^{\tilde K/2}$, is a vector of probability.
We bound from below the difference between $F_{\ga}(p^{(\nu)})$ and $F_{\ga}(p)$ in the next lemma, whose proof is postponed at the end of the section.
\begin{lem}\label{lem:separation}
For any $\ga \in (0,2)$, $\ga\neq 1$, and any $\nu \in \Vcal \{-1,1\}^{\tilde K/2}$, we have
\begin{equation*}
   |F_{\ga}(p^{(\nu)}) - F_{\ga}(p)|  \geq  C \sum_{k=1}^{\tilde K/2} p_{2k}^{\ga-2} \de_{2k}^2 = : R
\end{equation*}for a constant $C>0$ depending only on $\ga$.
\end{lem}

%\textcolor{red}{cette démo n'est pas très lisible pour le lecteur, car elle utilise un THM et des notations de \cite{rohde2020geometrizing} qui sont inconnues ici.} 
We will show that it is hard to know if the data come from $p$ or a uniform mixture of the $p^{(\nu)}$, $\nu\in \Vcal$. We do so by using Theorem A.1  of \cite{rohde2020geometrizing}, with the notations  of \cite{rohde2020geometrizing}. For any fixed $\al$-LDP interactive mechanism $Q$, we write  $Q^n := (Qp)^n \in \textup{conv}\pa{Q\Pcal^{(n)}_{\leq F_{\ga}(p)}}$ and $\overline Q^n := 2^{-\tilde{K}/2} \sum_{\nu \in \Vcal} (Qp^{(\nu)})^n \in \textup{conv}\pa{Q\Pcal^{(n)}_{\geq F_{\ga}(p) + R}}$. With the notations of \cite{rohde2020geometrizing} and standard relations between probability metrics, we have that the upper affinity satisfies 
\begin{equation}\label{quantitélourde}
    \eta^{(n)}_A(Q,R) \geq \pi(Q^n, \overline{Q}^n) = 
    1 - d_{TV}(Q^n, \overline{Q}^n) \geq 1 - \sqrt{D_{kl}(Q^n, \overline{Q}^n)/2} \enspace.
\end{equation}
%We bound the Kullback-Leibler divergence between $Q^n$ and $\overline{Q}^n$ in the next lemma.

%\begin{lem}\label{lem:informationTheoretic}
We can bound the KL-divergence $D_{kl} (Q^n, \overline Q^n)$ as in the proof of Theorem 4.2 in \cite{butucea2020interactive}, and have 
\begin{equation*}
   D_{kl} (Q^n, \overline Q^n) \leq  \frac{n (e^{2\al}-e^{-2\al})^2}{4}\|\de\|_2^2 \enspace.
\end{equation*}

Hence, it suffices to choose a $\de$ satisfying the condition
\begin{equation}\label{condLB:proof}
    \|\de\|_2^2 \leq \frac{2}{n (e^{2\al}-e^{-2\al})^2} \enspace,
\end{equation} to have $\eta^{(n)}_A(Q,R) \geq \frac{1}{2}$. Denoting $\Delta_A^{(n)}(Q,\eta) := \sup\{\Delta \geq 0\, :\, \eta_A^{(n)}(Q,\Delta) > \eta\}$ as in \cite{rohde2020geometrizing}, we will get for any $\eta \in (0,1/2)$,
\begin{equation*}
    \Delta_A^{(n)}(Q,\eta)  \geq  R
\end{equation*}where $R$ is defined in Lemma \ref{lem:separation} above.
It will then follow from Theorem A.1 of \cite{rohde2020geometrizing} that
\begin{equation*}
    \underset{Q}{\textup{inf}} \, \underset{\hat{F}}{\textup{inf}} \, \underset{p  \in \Pcal}{\textup{sup}} \, \E \cro{ ( \hat{F} - F_{\ga}(p))^2 }  \geq  \left(\frac{R}{2}\right)^2 \frac{\eta}{2} \enspace,
\end{equation*}for any $\eta \in (0,1/2)$. Taking $\eta = 1/4$ we will have 
\begin{equation*}
    \underset{Q}{\textup{inf}} \, \underset{\hat{F}}{\textup{inf}} \, \underset{p  \in \Pcal}{\textup{sup}} \, \E \cro{ ( \hat{F} - F_{\ga}(p))^2 }  \geq   \frac{C^2}{32} \pa{\sum_{k=1}^{\tilde K/2} p_{2k}^{\ga-2} \de_{2k}^2}^2 \enspace.
\end{equation*}
To choose a $\de$ fulfilling \eqref{condLB:proof}, we consider two cases according to the values of $K$.

\textit{$1^{\circ}.$ In the case where $K < n (e^{2\al}-e^{-2\al})^2$,} we choose $\tilde K = K$, and take $\de_k= (4\sqrt{K n } (e^{2\al}-e^{-2\al}))^{-1}$, $k\in[K]$. We take $p_k = 2\de_k$, $k\in[K-2]$, and the remaining $p_{K-1}, p_{K}  \geq 2\de_k$ so that $p$ is a vector of probability (i.e. $\sum_k p_k =1$). This gives
\begin{equation*}
    \underset{Q}{\textup{inf}} \, \underset{\hat{F}}{\textup{inf}} \, \underset{p  \in \Pcal}{\textup{sup}} \, \E \cro{ ( \hat{F} - F_{\ga}(p))^2 }  \geq  \frac{C^2}{32} \left(\frac{ 2^{\ga -2} \left[(K/2)-1\right] }{(4\sqrt{K n } (e^{2\al}-e^{-2\al}))^{\ga}}\right)^2 \geq \frac{C^2 2^{-2\ga}}{8192} \frac{  K^{2-\ga}}{( (e^{2\al}-e^{-2\al})^2 n   )^{\ga}}  \enspace,
\end{equation*}where we used $(K/2)-1 \geq K/4$ with $K\geq 4$. This corresponds to the right term of both lower bounds announced in Theorem 2.7.

\textit{$2^{\circ}.$ In the case where $K \geq n (e^{2\al}-e^{-2\al})^2$,} we separate our analysis in two ranges of values of $\ga$.\\
If $\ga \in(0,1)$, we take $\tilde K = K$, and  $\de_k= (2K)^{-1} $  and $p_k = 2\de_k$ for all $k\in[K]$. This leads to
\begin{equation*}
    \underset{Q}{\textup{inf}} \, \underset{\hat{F}}{\textup{inf}} \, \underset{p  \in \Pcal}{\textup{sup}} \, \E \cro{ ( \hat{F} - F_{\ga}(p))^2 }  \geq   \frac{C^2}{32} \left(\frac{ 2^{\ga -2} (K/2) }{(2 K)^{\ga}}\right)^2 \geq \frac{C^2}{2048}  K^{2(1-\ga)} \enspace,
\end{equation*}
which matches the first term of the  lower bound for $\ga \in(0,1)$  in the theorem.\\
If $\ga \in(1,2)$, let $\tilde K$ be the smallest even integer satisfying $\tilde K \geq n (e^{2\al}-e^{-2\al})^2$ and $\tilde K\geq 4$. We set $\de_k= (8\sqrt{ \tilde K n } (e^{2\al}-e^{-2\al}))^{-1}$ for $k\in[\tilde K]$. We  choose $p_k = 2\de_k$ for $k\in[\tilde K-2]$, and  $p_k \geq 2 \de_k$ for $k\geq \tilde K-1$ such that $p$ is a vector of probability. Then
\begin{align*}
    \underset{Q}{\textup{inf}} \, \underset{\hat{F}}{\textup{inf}} \, \underset{p  \in \Pcal}{\textup{sup}} \, \E \cro{ ( \hat{F} - F_{\ga}(p))^2 }  &\geq  \frac{C^2}{32} \left(\frac{ 2^{\ga -2} \left[(\tilde K/2)-1\right] }{(8\sqrt{ \tilde K n } (e^{2\al}-e^{-2\al}))^{\ga}}\right)^2\\
    & \geq \frac{C^2 4^{-2\ga}}{8192} \frac{  \tilde K^{2-\ga}}{( (e^{2\al}-e^{-2\al})^2 n   )^{\ga}}  \\
    &\geq \frac{C^2 4^{-2\ga}}{8192} ( (e^{2\al}-e^{-2\al})^2 n   )^{2(1-\ga )} \enspace, 
\end{align*}which corresponds to the first term of the lower bound for $\ga \in (1,2)$ in the theorem.

The proof of Theorem 2.7 is complete.\hfill $\square$

%\textcolor{brown}{On doit prendre $K \geq n \alpha^2$. On peut choisir de perturber $K' = n \alpha^2$ valeurs et obtenir la borne $\frac{\widetilde C'}{(\al^2 n)^{2(\ga-1)}}$, mais on ne peut avoir $\frac{K^{2-\ga}}{(\al^2n)^{\ga}}$ autrement que pour $K=n \alpha^2$...} \ya{oui}

%%%%%%%%%%%%%%%%%%%%%%%%%%%%%%%%%%
%%%%%%%%%%%%%%%%%%%%%%%%%%%%%%%%%%%

%\subsection{Proof of lemmas for lower bounds}

\bigskip

\noindent \textbf{Proof of Lemma \ref{lem:separation}.} We have
\begin{equation}\label{eq:proof:sep}
    F_{\ga}(p^{(\nu)}) - F_{\ga}(p) = \sum_{k=1}^{\tilde K/2} \Big{[}(p_{2k}+ \nu_k \de_{2k})^{\ga} + (p_{2k}-\nu_k\de_{2k})^{\ga} - 2 p_{2k}^{\ga} \Big{]}
\end{equation}
Denoting $f(x)=x^{\ga}$ and using Taylor's formula, we have for any real $Y >0$,
\begin{align*}
f(Y) &= f(p_{2k}) + f'(p_{2k}) (Y-p_{2k}) + f''(w_Y)\frac{(Y-p_{2k})^2}{2}
\end{align*}
where $w_Y$ lies between $Y$ and $p_{2k}$.
We take $Y= p_{2k}+ \nu_k \de_{2k}$ and $\widetilde Y =p_{2k} - \nu_k \de_{2k}$ to get
\begin{align}\label{proof:lem:nonpositive}
    f(Y) + f(\widetilde Y) - 2p_{2k}^{\ga} &= f''(w_Y)\frac{(Y-p_{2k})^2}{2} + f''( w_{\widetilde Y})\frac{(Y-p_{2k})^2}{2} \nonumber \\
    &=  \ga(\ga-1) (w_{Y}^{\ga-2} + w_{\widetilde Y}^{\ga-2}) \frac{\de_{2k}^2}{2}
\end{align}
Since $w_{Y} \vee w_{\widetilde Y} \leq p_{2k} + \de_{2k}$ with  $0\leq \de_{2k} \leq p_{2k}/2$, and $\ga \in(0,2)$, we have
\begin{equation*}
   w_{Y}^{\ga-2} \wedge w_{\widetilde Y}^{\ga-2} \geq   (p_{2k}+ \de_{2k})^{\ga-2}  \geq (2p_{2k})^{\ga-2} \enspace .
\end{equation*}
Hence, for $\ga \in (1,2)$  ,
\begin{equation*}
    f(Y) + f(\widetilde Y) - 2p_{2k}^{\ga} \geq \ga(\ga-1)  (2p_{2k})^{\ga-2} \de_{2k}^2
\end{equation*}which leads to the desired lower bound of \eqref{eq:proof:sep}.  For $\ga \in(0,1)$, we deduce from \eqref{proof:lem:nonpositive} that all terms of the sum \eqref{eq:proof:sep} are non-positive and satisfy
\begin{equation*}
    f(Y) + f(\widetilde Y) - 2p_{2k}^{\ga} \leq \ga(\ga-1)  (2p_{2k})^{\ga-2} \de_{2k}^2\enspace.
\end{equation*}So, the absolute value of the sum \eqref{eq:proof:sep} can be lower bounded as announced in the lemma. \hfill $\square$

\bibliographystyle{plain}

\end{document}